\documentclass[12pt]{article}

\usepackage[mathscr]{eucal}
\usepackage{amsmath, amsthm, amssymb, latexsym, graphicx,courier,makeidx,subeqnarray,color}
\usepackage{bm}
\usepackage{color}

\newtheorem{Theorem}{Theorem}
\newtheorem{Definition}[Theorem]{Definition}
\newtheorem{Lemma}[Theorem]{Lemma}

\newtheorem{Problem}[Theorem]{Problem}

\newtheorem{Corollary}{Corollary}[section]
\newtheorem{Remark}{Remark}[section]

\def\real{\mathbb{R}}
\def\bar{\overline}

\newcommand{\bu}{\mbox{\boldmath{$u$}}}
\newcommand{\bv}{\mbox{\boldmath{$v$}}}
\newcommand{\bw}{\mbox{\boldmath{$w$}}}
\newcommand{\bx}{\mbox{\boldmath{$x$}}}

\newcommand{\fb}{\mbox{\boldmath{$f$}}}

\newcommand{\hb}{\mbox{\boldmath{$h$}}}
\newcommand{\bxi}{\mbox{\boldmath{$\xi$}}}
\newcommand{\beeta}{\mbox{\boldmath{$\eta$}}}

\newcommand{\bsigma}{\mbox{\boldmath{$\sigma$}}}
\newcommand{\btau}{\mbox{\boldmath{$\tau$}}}
\newcommand{\bvarepsilon}{\mbox{\boldmath{$\varepsilon$}}}
\newcommand{\bnu}{\mbox{\boldmath{$\nu$}}}
\newcommand{\bzeta}{\mbox{\boldmath{$\zeta$}}}

\newcommand{\bzero}{\mbox{\boldmath{$0$}}}

\newcommand{\bz}{\mbox{\boldmath{$z$}}}

\def\lista#1
{{ \itemindent 0.0cm \labelsep .2cm \leftmargin 0.8cm \rightmargin
		0.0cm \labelwidth 0.6cm \topsep 0.0mm
		\parsep 0.0mm
		\itemsep 0.0mm
		\begin{list}{}
			{ \setlength{\leftmargin}{.8cm} \setlength{\rightmargin}{0.0cm}
				\setlength{\parsep}{0.0mm} \setlength{\topsep}{.0mm}
				\setlength{\parskip}{.0cm} \setlength{\itemsep}{.0cm} }
			{#1}\end{list}} }

\begin{document}
	
\title{
Well-posedness of Constrained Evolutionary Differential 
Variational--Hemivariational Inequalities \\
with Applications
	\thanks{\, 
	The project has received funding from the European Union's Horizon 2020 Research and Innovation Programme under the Marie Sk{\l}odowska-Curie grant agreement No. 823731 CONMECH. 
	It is supported by 
	Natural Science Foundation of Guangxi (Grant No: 2018GXNSFAA281353), 
	Beibu Gulf University Project No. 2018KYQD06, and the projects financed by the Ministry of Science and
	Higher Education of Republic of Poland under
	Grants Nos. 4004/GGPJII/H2020/2018/0 and
	440328/PnH2/2019, and the National Science Centre of Poland under Project No. 2021/41/B/ST1/01636.
}}

\author{
	Stanis{\l}aw Mig\'orski 
	\footnote{\, College of Applied Mathematics, Chengdu University of Information Technology, Chengdu 610225, Sichuan Province, P.R. China, and Jagiellonian University in Krakow, Chair of Optimization and Control, ul. Lojasiewicza 6, 30348 Krakow, Poland.
	Tel.: +48-12-6646666. E-mail address: stanislaw.migorski@uj.edu.pl.}
}

\renewcommand{\thefootnote}{\fnsymbol{footnote}}

\date{}
\maketitle
\thispagestyle{empty}

\begin{abstract}
\noindent
A system of a first order history-dependent evolutionary  variational-hemivariational inequality with unilateral constraints coupled with a nonlinear ordinary differential equation in a Banach space is studied. 
Based on 
a fixed point theorem for history dependent operators, 
results on the well-posedness of the system are proved. 
Existence, uniqueness, continuous dependence of the solution
on the data, and the solution regularity are established.
Two applications of dynamic problems from contact mechanics illustrate the abstract results. 
First application is a unilateral viscoplastic frictionless contact problem which leads to a hemivariational inequality
for the velocity field, 
and the second one deals with a viscoelastic frictional
contact problem which is described by a variational inequality.


\end{abstract}








\section{Introduction}\label{Intro}
In this paper we study the following system of an  evolutionary variational-hemivariational inequality coupled with an ordinary differential equation in a Banach space.

Find ${\mathtt x} \colon (0, T) \to E$ 
and $w \colon (0, T) \to V$ with 
$w(t) \in K$ for a.e. $t \in (0, T)$ 
such that
\begin{eqnarray}
&&\hspace{-0.5cm}
{\mathtt x}'(t) = F(t,{\mathtt x}(t), w(t), (Sw)(t))
\ \ \mbox{\rm for a.e.} \ \ t \in (0,T), 
\label{001} 
\\[1mm]
&&\hspace{-0.5cm}
\langle w'(t) + A(t, {\mathtt x}(t), w(t)) 
+ (R_1 w)(t), v - w(t) \rangle 
+ \, j^0 (t, {\mathtt x}(t), 
(R_2w)(t), Mw(t); Mv - Mw(t)) \nonumber 
\\[1mm]
&&
\quad \ \, 
+ \, \varphi (t, {\mathtt x}(t), (R_3w)(t), Mv) 
- \varphi(t, {\mathtt x}(t), (R_3w)(t), Mw(t)) \ge 
\langle f(t, {\mathtt x}(t)), v - w(t) \rangle 
\nonumber
\\[1mm]
&&\qquad 
\ \ \ \ \mbox{\rm for all} \ v \in K, \ \mbox{\rm a.e.} \ t \in (0,T), 
\label{002} \\
&&\hspace{-0.5cm}
{\mathtt x}(0) = {\mathtt x}_0, \ \ w(0) = w_0 .
\label{003}
\end{eqnarray}
Here, $E$ is a Banach space, 
and $(V, H, V^*)$ is an evolution triple of spaces. 
Problem (\ref{001})--(\ref{003}) represents a system which couples the ordinary differential equation (\ref{001}) 
with the constrained variational-hemivariational
inequality (\ref{002}), supplemented with the 
initial conditions (\ref{003}). 
Following the terminology in~\cite{Anh,LMZ2017,Liu1,MZJOGO2018}, 
we refer to (\ref{001})--(\ref{003}) as an evolutionary differential variational-hemivariational inequality. 
%
The inequality (\ref{002}) couples with equation (\ref{001}) through the solution 
$w$ and its history value $Sw$, $j^0$ denotes the generalized directional 
derivative in the last variable of a nonconvex and nondifferentiable function 
$j \colon (0, T) \times E \times Z \times X \to \real$, 
$\varphi \colon (0, T) \times E \times Y 
\times X \to \real$ is a function convex in the last variable, $M \colon V\rightarrow X$ is an affine linear  operator, 
$S$, $R_1$, $R_2$, $R_3$ are the so-called
history-dependent operators, 
and $K \subset V$ represents a set of constraints. 
%
For a history-dependent operator, the current value 
for a given function at the time instant $t$ depends 
on the va\-lues of the function at the time instants 
from $0$ to $t$. For this reason, the probem 
(\ref{001})--(\ref{003}) is also called a history-dependent differential variational-hemivariational inequality.
The main features of the system are the presence of the constraints set $K$ 
and the strong coupling between the differential equation and the inequality which leads to a complex time dependent dynamics. 
Problem (\ref{001})--(\ref{003}) is new and, to the best of our knowledge, has not been studied in the literature. 
Our main existence and uniqueness result, 
in its particular case,  
solves an open problem stated  in~\cite[Section~10.4]{SOFMIG}, 
and extends the recent result in~\cite[Theorem~20]{MigorskiBiao} 
obtained for a purely variational-hemivariational inequality.

There are two main goals of the paper. 
First, we establish a well-posedness result which includes existence, uniqueness and continuous dependence 
of solution to the system (\ref{001})--(\ref{003})
on the data. 
In particular, 
if $f$ is assumed to be independent of ${\mathtt x}$, 
then we show, see~Theorem~\ref{Theorem1}, that 
the map 
$({\mathtt x}_0, w_0, f) \mapsto ({\mathtt x}, w)$
is Lipschitz continuous. 
The second goal of the paper is to provide two applications to new classes of dynamic problems in contact mechanics.
In the first application we study a dynamic unilateral viscoplastic frictionless contact problem involving  
a nonmonotone Clarke subdifferential boundary condition.
Its
weak formulation is a hemivariational inequality coupled with an ordinary differential equation.
The second application concerns a dynamic frictional viscoelastic contact problem with adhesion which leads 
to a variational inequality combined with a differential equation on the contact surface.


Only special versions of the system (\ref{001})-(\ref{003}) 
have been explored in the literature.  
%
%
For instance, the evolutionary variational-hemivariational inequalities without a constraint set $K$, 
and with and without history-dependent operators and their variants are discussed in~\cite[Chapter~7]{SOFMIG}. 
If $K=V$, $\varphi = 0$, and 
$A$, $j$, $f$ are independent of ${\mathtt x}$,  
then (\ref{002})-(\ref{003}) reduces to the evolution hemivariational inequality considered 
in~\cite{MOgorzaly}. 
If the operator $A$ and a locally Lipschitz function $j$ are assumed to be independent of 
${\mathtt{x}}$ and $R_1$, then the problem reduces to the parabolic hemivariational inequality studied, 
for example, in~\cite{Liu2008, Miettinen,MO2004}. 
Very recently the system has been treated in~\cite{ZM2021}  
with $K=V$, $\varphi = 0$, and $F$, $j$ independent of the history-dependent operators, and $f$ independent of 
the variable ${\mathtt x}$.
A general dynamic variational-hemivariational
inequality with history-dependent operators 
without constraints can be found in~\cite{HMS2017}.
All aforementioned papers treat the special case of (\ref{002}) with $K = V$. 
In contact mechanics, the ordinary differential 
equation (\ref{001}) supplementing (\ref{002}) 
appears na\-tu\-rally. 
For example, the system is met in
rate-type viscoplastic constitutive laws, 
including the internal state variables 
in viscoplasticity,
see~\cite[Chapter~3]{SOFMIG}, 
and in modeling of additional effects 
like wear and adhesion phenomena in contact problems, 
see~\cite[Chapter~11]{SST}, \cite[Chapter~2]{SHS}
and the references therein.

We remark that the notion of a history-dependent operator was introduced in~\cite{SM1} and used in a several recent  papers~\cite{HMS2017,Migorski2020,MOS13,MOS18,MOgorzaly,
	SHM,SM1,SMH2018,SP,SX}.
Differential variational inequalities was initiated and first systematically discussed in~\cite{Pang} in finite dimensional spaces.
Note also that a related result on a class of differential hemivariational inequalities which consists of a hemivariational inequality of parabolic type combined 
with a nonlinear evolution equation has been delivered 
by the Rothe method in~\cite{MZJOGO2018}.
The differential parabolic-parabolic variational
inequalities are examined in~\cite{Anh} by using the technique of measure of noncompactness.
Results on partial differential variational inequalities with nonlocal boundary conditions can be found in~\cite{LMZ2017}. 
The literature on hemivariational inequalities has been significantly enlarged in the last forty years, 
see monographs~\cite{CLM,Go11.I,MOSBOOK, NP, SOFMIG}, 
for analysis of various classes of such inequalities,
see e.g.~\cite{Bartosz,Gwinner,HMS2017,KULIG,
Liu2008,LMZ2017,Liu1,MMM,MHZ2020,Migorski2021,MZJOGO2018,ZM2021}. 
By our approach, we relax the assumptions on similar problems treated in the aforementioned papers and 
considerably improve some results by allowing 
the history-dependent operators to appear 
in the convex term and in the generalized directional derivative of a nonconvex potential.

%

The paper is structured as follows. 
In Section~\ref{Prelim}, we recall some prerequisites needed throughout this paper. 
In Section~\ref{s1}, we state and prove the main results on the well-posedness of the system (\ref{001})--(\ref{003}). 
Next, in Section~\ref{Application1}, we consider 
a model of contact which leads to a differential hemivariational inequality for the velocity field 
of the form (\ref{001})--(\ref{003}). 
Finally, for the second contact problem in Section~\ref{Application2}, 
under appropriate hypotheses on the data, we prove the well-posedness for a differential variational inequality with constraints which is a weak form of the contact model.


\section{Preliminaries}\label{Prelim}

\noindent
%
Let $(X, \| \cdot \|_X)$ be a normed space, 
$X^*$ be the dual of $X$ and 
$\langle\cdot,\cdot\rangle_{X^*\times X}$
denote the duality brackets for the pair $(X^*, X)$. 
For simplicity, when no confusion arises, we often skip the subscripts. 
For Banach spaces $X$, $Y$, we will use the notation 
${\mathcal L}(X, Y)$ for the set of all linear continuous
operators from $X$ to $Y$.
The notation $\| A \|$ stands for the operator norm 
of $A \colon X \to Y$ in ${\mathcal L(X, Y)}$.
%
Given a set $D \subset X$, we write  
$\| D \|_X = \sup \{ \| x \|_X \mid x \in D \}$.
We denote by $C([0, T]; X)$ the space of continuous
functions on $[0, T]$ with values in $X$.
It is well known that if $X$ is a Banach space, 
then $C([0, T]; X)$ is also a Banach space.

Let $h \colon X \to \real$ be a locally Lipschitz function. The
{\rm generalized (Clarke) directional derivative} of $h$ at the point
$x \in X$ in the direction $v \in X$ is defined
by
\begin{equation*}
h^{0}(x; v) = \limsup_{y \to x, \ \lambda \downarrow 0}
\frac{h(y + \lambda v) - h(y)}{\lambda}.
\end{equation*}
The {\rm generalized subgradient} of $h$ at $x$
is a subset of the dual space $X^*$ given by
\begin{equation*}
\partial h (x) = \{\, \zeta \in X^* \mid h^{0}(x; v) \ge
{\langle \zeta, v \rangle} \ \mbox{for all} \ v \in X \, \}.
\end{equation*}
A locally Lipschitz function $h$ is said to be {\rm regular} 
(in the sense of Clarke) at the point $x \in X$ if for all $v \in X$ the 
derivative $h' (x; v)$ exists and 
$h^0(x; v) = h'(x; v)$.
Throughout the paper, the generalized derivative and the generalized subgradient are always taken
with respect to the last variable of a given function.

We say that a map $A \colon X \to X^*$ 
is demicontinuous if the function 
$u \mapsto \langle A u, v \rangle_{X^* \times X}$
is continuous for all $v \in X$, i.e., $A$ is continuous as a mapping from $X$ to $X^*$ endowed with $w^*$-topology.
It is hemicontinuous, if for all $u$, $v$, 
$w \in X$, the function
$t \mapsto \langle A(u+tv), w \rangle_{X^* \times X}$ is continuous on $[0,1]$.
A map $A$ is monotone if $\langle Au - Av, u - v 
\rangle_{X^* \times X} \ge 0$ for all $u$, $v \in X$. 
It is strongly monotone with constant $c > 0$ 
if $\langle Au - Av, u - v \rangle_{X^* \times X} 
\ge c \, \|u-v\|^2_X$ for all $u$, $v \in X$. 
For more details, we refer to~\cite{Clarke,DMP1,DMP2,MOSBOOK}.

The following definition and a fixed point result 
(being a consequence of the Banach contraction 
principle) can be found 
in~\cite[Definition~30 and Theorem~67]{SOFMIG}, 
respectively.

\begin{Definition}\label{historyop}
Let	${\mathbb X}$ and ${\mathbb Y}$ be normed spaces. 
An operator 
${\mathbb S} \colon L^2(0,T; {\mathbb X}) \to
L^2(0,T; {\mathbb Y})$ is called a history-dependent operator if there exists $M > 0$ such that
$$
\| ({\mathbb S} u_1)(t)- ({\mathbb S}u_2)(t) 
\|_{{\mathbb Y}}
\le M \int_0^t \| u_1(s) -u_2(s)\|_{{\mathbb X}} \, ds
$$
for all $u_1$, $u_2 \in L^2(0,T; {\mathbb X})$, 
a.e. $t \in (0, T)$. 
\end{Definition} 	
%
%
\begin{Lemma}\label{CONTR}
	Let ${\mathbb X}$ be a Banach space and $0 < T < \infty$.
	Let $F \colon L^2(0,T; {\mathbb X}) \to L^2(0,T; {\mathbb X})$ be an operator such that
	\begin{equation*}
	\| (F \eta_1)(t) - (F \eta_2)(t)\|^2_{\mathbb X} \le c
	\int_0^t \| \eta_1(s) - \eta_2(s)\|^2_{\mathbb X} \, ds
	\end{equation*}
	for all \, $\eta_1$, $\eta_2 \in L^2(0,T;{\mathbb X})$, 
	a.e.\ $t \in (0,T)$ with a constant $c > 0$. Then $F$ has a unique fixed point in 
	$L^2(0, T; {\mathbb X})$, 
	i.e., there exists a unique
	$\eta^* \in  L^2(0,T;{\mathbb X})$ such that 
	$F \eta^* = \eta^*$.
\end{Lemma}

\section{Well-posedness of the abstract problem}\label{s1}

In this section we provide results on the well-posedness
of the differential variational-hemiva\-ria\-tio\-nal inequalities.
We shall establish existence and uniqueness of solution 
and prove continuous dependence of solution on 
the initial conditions and the function $f$. 

Let $E$, $U$, $X$, $Y$ and $Z$ be Banach spaces 
and 
$(V, H, V^*)$ be an evolution triple of spaces. 
The latter means that $V$ is a separable reflexive Banach space and $H$ is a separable Hilbert space such that the embedding $V \subset H$ is continuous and dense.
Then $H$ is embedded continuously and densely 
in $V^*$, and the duality brackets 
$\langle \cdot, \cdot \rangle$ for the pair 
$(V^*, V)$ and the inner product $\langle \cdot, \cdot \rangle_H$ on $H$ coincide on $H \times V$. 
In what follows we denote by $\| \cdot \|$ the norm in $V$.
We set 
$$
{\mathbb{W}} = \{ \, w \in L^2(0, T; V) \mid 
w' \in L^2(0,T; V^*) \, \},
$$
where the time derivative $w'$ is understood 
in the distributional sense. 
It is known that
$L^2(0,T;V)^* \simeq L^2(0,T; V^*)$ 
and the space ${\mathbb{W}}$ endowed with the norm
$\| w \|_{{\mathbb{W}}} = 
\| w \|_{L^2(0, T; V)} + \| w' \|_{L^2(0, T; V^*)}$
is a separable reflexive Banach space. 

We deal with the following system.

\begin{Problem}\label{DVHI} 
Find $ {\mathtt x}\in H^1(0,T; E)$ and 
$w \in {\mathbb{W}}$
such that $w(t) \in K$ for a.e. $t \in (0, T)$ 
and
\begin{equation*}
\begin{cases}
\displaystyle
{\mathtt x}'(t) = F(t,{\mathtt x}(t), w(t), (Sw)(t))
\ \ \mbox{\rm for a.e.} \ \, t \in (0,T), \\[1mm]
\langle w'(t) + A(t, {\mathtt x}(t), w(t)) 
+ (R_1 w)(t), v - w(t) \rangle 
+ 
\, j^0 (t, {\mathtt x}(t), 
(R_2w)(t), Mw(t); Mv - Mw(t)) \\[1mm]
\ \ \quad \qquad 
+ \, \varphi (t, {\mathtt x}(t), (R_3w)(t), Mv) 
- \varphi(t, {\mathtt x}(t), (R_3w)(t), Mw(t)) \ge 
\langle f(t, {\mathtt x}(t)), v - w(t) \rangle 
\\[1mm]
\qquad\qquad\qquad \qquad
\ \ \mbox{\rm for all} \ v \in K, \ \mbox{\rm a.e.} \ t \in (0,T), \\
{\mathtt x}(0) = {\mathtt x}_0, \ \ w(0) = w_0 .
\end{cases}
\end{equation*}
\end{Problem}
We need the following hypotheses on the data of 
Problem~\ref{DVHI}.

\smallskip

\noindent 
$\underline{H({F})}:$ \quad 
$\displaystyle 
F \colon (0, T) \times E \times V \times U \to E$ 
is such that 
\smallskip

\lista{
\item[(a)]  
$F(\cdot, {\mathtt x}, v, u)$ belongs to $L^2(0,T; E)$ 
for all ${\mathtt x} \in E$, $v \in V$, $u\in U$. 
\smallskip
\item[(b)]
$\| F(t, {\mathtt x}_1, v_1, u_1) - 
F(t, {\mathtt x}_2, v_2, u_2)\|_E
\le L_F \, (\| {\mathtt x}_1-{\mathtt x}_2\|_E 
+ \| v_1-v_2 \| + \| u_1 - u_2\|_U)$ 
for all ${\mathtt x}_1$, ${\mathtt x}_2 \in E$, 
$v_1$, $v_2 \in V$, $u_1$, $u_2 \in U$, 
a.e. $t \in (0, T)$. 
}

\smallskip

\smallskip
\noindent 
$\underline{H({A})}:$ \quad 
$\displaystyle 
A \colon (0, T) \times E \times V \to V^*$ is such that 
\smallskip

\lista{
\item[(a)]  
$A(\cdot, {\mathtt x}, v)$ is measurable for all 
${\mathtt x} \in E$, $v \in V$. 
\smallskip
\item[(b)]
$A(t,\cdot, v)$ is continuous for all $v \in V$, 
a.e. $t \in (0, T)$. 
\smallskip	
\item[(c)]
$\| A(t, {\mathtt x}, v) \|_{V^*} 
\le a_0(t) + a_1 \| {\mathtt x} \|_E 
+ a_2 \| v \|$ for all ${\mathtt x} \in E$, 
$v \in V$, a.e. 
$t \in (0, T)$ with $a_0 \in L^2(0,T)$, $a_0$, $a_1$, $a_2 \ge 0$. 
\smallskip
\item[(d)] 
$A(t, {\mathtt x}, \cdot)$ is demicontinuous and there are constants $m_A >0$, 
${\bar{m}}_A \ge 0$ such that
$$\langle A (t, {\mathtt x}_1, v_1) 
- A(t, {\mathtt x}_2, v_2), v_1 - v_2 
\rangle
\ge m_A \, \| v_1 - v_2 \|^2
- {\bar{m}}_A \, \|{\mathtt x}_1 - {\mathtt x}_2 \|_E 
\, \| v_1 - v_2\|
$$
for all ${\mathtt x}_1$, ${\mathtt x}_2 \in E$, 
$v_1$, $v_2 \in V$, a.e. $t \in (0, T)$. 
}

\smallskip

\smallskip

\noindent 
$\underline{H(j)}:$ \quad 
$j \colon (0, T) \times E \times Z \times X \to \real$ is such that 

\smallskip

\lista{
\item[(a)]  
$j(\cdot, {\mathtt x}, z, v)$ is measurable for all 
${\mathtt x} \in E$, $z \in Z$, $v \in X$. 
\smallskip
\item[(b)]
$j(t, \cdot, \cdot, v)$ is continuous for all 
$v \in X$, a.e. $t \in (0, T)$. 
\smallskip
\item[(c)]
$j(t, {\mathtt x}, z, \cdot)$ 
is locally Lipschitz for all 
${\mathtt x} \in E$, $z \in Z$, a.e. $t \in (0, T)$. 
\smallskip
\item[(d)]
$\| \partial j(t, {\mathtt x}, z, v) \|_{X^*} \le 
c_{0j} (t) + c_{1j} \| {\mathtt x} \|_E + 
c_{2j} \| z \|_Z + c_{3j} \| v \|_X$ \smallskip
for all ${\mathtt x} \in E$, $z \in Z$, 
$v \in X$, a.e. $t \in (0, T)$ 
with $c_{0j} \in L^2(0, T)$, 
$c_{0j}$, $c_{1j}$, $c_{2j}$, $c_{3j} \ge 0$. 
\smallskip
\item[(e)] 
$j^0(t, {\mathtt x}_1, z_1, v_1; v_2 - v_1) 
+ j^0(t, {\mathtt x}_2, z_2, v_2; v_1 - v_2) \\[1mm]
~~ \qquad \qquad \qquad \qquad 
\le m_j \, \| v_1 - v_2 \|^2_X + 
{\bar{m}}_{j}
\, (\| {\mathtt x}_1 - {\mathtt x}_2\|_E 
+ \| z_1 - z_2 \|_Z ) \, \| v_1 - v_2 \|_X$ 
\\[2mm]
for all ${\mathtt x}_i\in E$, $z_i \in Z$, 
$v_i \in X$, 
$i = 1$, $2$, a.e.\ $t \in (0, T)$ with $m_j$, ${\bar{m}}_j \ge 0$.
}

\smallskip

\noindent 
$\underline{H(\varphi)}:$ \quad 
$\varphi \colon (0, T) \times E \times Y 
\times X \to \real$ is such that 

\smallskip

\lista{
\item[(a)]  
$\varphi(\cdot, {\mathtt x}, y, v)$ 
is measurable for all 
${\mathtt x} \in E$, $y \in Y$, $v \in X$. 
\smallskip
\item[(b)]
$\varphi(t, \cdot, \cdot, v)$ is continuous for all 
$v \in X$, a.e. $t \in (0, T)$. 
\smallskip
\item[(c)]
$\varphi(t, {\mathtt x}, y, \cdot)$ is convex and lower semicontinuous for all ${\mathtt x} \in E$, $y \in Y$, a.e. 
$t \in (0, T)$. \smallskip
\item[(d)]
$\| \partial \varphi(t, {\mathtt x}, y, v) \|_{X^*} \le  c_{0\varphi} (t) +
c_{1\varphi} \| {\mathtt x} \|_E 
+ c_{2\varphi} \| y \|_Y 
+ c_{3\varphi} \| v \|_X$ 
for all ${\mathtt x} \in E$, $y \in Y$, $v \in X$, a.e. 
$t \in (0, T)$ with $c_{0\varphi} \in L^2(0, T)$, 
$c_{0\varphi}$, $c_{1\varphi}$, $c_{2\varphi}$, $c_{3\varphi} \ge 0$. 
\smallskip
\item[(e)] 
$
\varphi (t, {\mathtt x}_1, y_1, v_2) - 
\varphi (t, {\mathtt x}_1, y_1, v_1) + 
\varphi (t, {\mathtt x}_2, y_2, v_1) - 
\varphi (t, {\mathtt x}_2, y_2, v_2) \\[1mm]
~~ \qquad \qquad \qquad \qquad
\le m_{\varphi} \, 
(\|{\mathtt x}_1 - {\mathtt x}_2 \|_E + 
\| y_1 - y_2 \|_Y ) \, \| v_1 - v_2 \|_X
$ \\[2mm]
for all ${\mathtt x}_i \in E$, $y_i \in Y$, 
$v_i \in X$, $i = 1$, $2$, 
a.e.\ $t \in (0, T)$ with 
$m_{\varphi} \ge 0$.
}

\smallskip

\smallskip

\noindent 
$\underline{H(K)}:$ \quad
$K$ is a closed and convex subset of $V$ with 
$0 \in K$.

\smallskip

\smallskip

\noindent
$\underline{H(M)}:$ \quad
$M \colon V\rightarrow X$ is such that

\lista{
\item[(a)]
$M$ is an affine operator.
\smallskip
\item[(b)]
the Nemitsky operator ${\mathcal M} \colon L^2(0,T;V) 
\to L^2(0,T;X)$ corresponding to $M$ is compact.
}

\smallskip

\smallskip

\noindent 
$\underline{H(f)}:$ \quad
$f \colon (0, T) \times E \to V^*$ is such that

\smallskip

\lista{
	\item[(a)]  
	$f(\cdot, {\mathtt x})$ belongs to $L^2(0, T; V^*)$ 
	for all ${\mathtt x} \in E$. 
	\smallskip 
	\item[(b)]
	$\| f(t, {\mathtt x}_1) - f(t, {\mathtt x}_2) \|_{V^*} \le L_f \, 
	\| {\mathtt x}_1-{\mathtt x}_2\|_E$ 
	for all ${\mathtt x}_1$, 
	${\mathtt x}_2 \in E$, a.e. $t \in (0, T)$ 
	with $L_f > 0$. 
}

\smallskip

\smallskip

\noindent 
$\underline{(H_1)}:$ \quad 
$m_j \| M_1 \|^2 < m_A$, where
$M_1 \colon V \to X$ 
is defined by $M_1 v = M v - M0$ for $v \in V$.

\smallskip

\smallskip

\noindent 
$\underline{(H_2)}:$ \quad 
${\mathtt x}_0 \in E$, $w_0 \in V$.

\smallskip

\smallskip

\noindent 
$\underline{(H_3)}:$ \quad 
$S \colon L^2(0,T;V) \to L^2(0,T;U)$,
${R}_1 \colon L^2(0, T; V) \to L^2(0, T; V^*)$, 

\smallskip

\qquad
${R}_2 \colon L^2(0, T; V) \to L^2(0, T; Z)$, 
${R}_3 \colon L^2(0, T; V) \to L^2(0, T; Y)$ 
are such that

\smallskip

\lista{
\item[(a)] 
$\displaystyle 
\| (S v_1)(t) - (S v_2)(t) \|_U \le 
c_{S} \int_0^t \| v_1(s) - v_2(s) \| \, ds$,
\\ 
\item[(b)] 
$\displaystyle 
\| ({R}_1 v_1)(t) - ({R}_1 v_2)(t) \|_{V^*} \le 
c_{R_1} \int_0^t \| v_1(s) - v_2(s) \| \, ds$,
\\ 
\item[(c)] 
$\displaystyle
\| ({R}_2 v_1)(t) - ({R}_2 v_2)(t) \|_{Z} \le 
c_{R_2} \int_0^t \| v_1(s) - v_2(s) \| \,ds$, \\ 
\item[(d)] 
$\displaystyle 
\| ({R}_3 v_1)(t) - ({R}_3 v_2)(t) \|_Y \le 
c_{R_3} \int_0^t \| v_1(s) - v_2(s) \| \, ds$
}
\quad for all $v_1$, $v_2 \in L^2(0, T; V)$, 
a.e.\ $t\in (0, T)$ with $c_S$, $c_{R_1}$, 
$c_{R_2}$, $c_{R_3} > 0$.

\begin{Remark}\label{REM1}
{\rm 
In hypothesis $H(M)$, the operator $M$
is an affine operator if and only if the operator
$M_1$ defined in $(H_1)$ is linear. 
The operator $M_1$ is called the linear part of $M$.
In $(H_1)$ and in what follows, the constant 
$\| M_1 \|$ denotes the norm in ${\mathcal L}(V, X)$ 
of the linear part $M_1$ of $M$.
It is clear that if $M \in {\mathcal L}(V, X)$, 
then $M_1 = M$.
}
\end{Remark}

\begin{Remark}\label{REM2}
{\rm
The following conditions are useful while checking 
the hypothesis $H(A)${\rm (d)}.
If $A(t, {\mathtt x}, \cdot)$
is strongly monotone with $m_A>0$ 
and $A(t,\cdot, v)$ is 
Lipschitz with $L_A>0$, i.e., 
\begin{eqnarray*}
&&
\langle A (t, {\mathtt x}, v_1) 
- A(t, {\mathtt x}, v_2), v_1 - v_2 \rangle
\ge m_A \| v_1 - v_2 \|^2  
\ \ \mbox{for all} \ \  {\mathtt x} \in E, 
\ v_1, v_2 \in V
, \\[2mm]
&&
\| A (t, {\mathtt x}_1, v) 
- A(t, {\mathtt x}_2, v) \|_{V^*} \le L_A\, 
\| {\mathtt x}_1 - {\mathtt x}_2\|_E  
\ \ \mbox{for all} \ \ 
{\mathtt x}_1, {\mathtt x}_2 \in E,
\ v \in V, 
\end{eqnarray*}
for a.e. $t \in (0,T)$,
then
\begin{equation*}
\langle A (t, {\mathtt x}_1, v_1) 
- A(t, {\mathtt x}_2, v_2), v_1 - v_2 \rangle
\ge m_A \, \| v_1 - v_2 \|^2
- L_A \, \|{\mathtt x}_1 - {\mathtt x}_2 \|_E
\| v_1 - v_2\|
\end{equation*}
for all ${\mathtt x}_1$, ${\mathtt x}_2 \in E$, 
$v_1$, $v_2 \in V$, a.e. $t \in (0, T)$. 
}
\end{Remark}

The main result on the well-posedness of 
Problem~$\ref{DVHI}$ reads as follows.

\begin{Theorem}\label{Theorem1}
Under hypotheses $H(F)$, $H(A)$, $H(j)$, $H(\varphi)$, $H(K)$, $H(M)$, $H(f)$, $(H_1)$--$(H_3)$, 
for each $({\mathtt x}_0, w_0) \in E \times V$,
Problem~$\ref{DVHI}$ has a unique solution
$({\mathtt x}, w) \in H^1(0,T; E) \times {\mathbb{W}}$ 
with $w(t) \in K$ for a.e. $t \in (0, T)$.
Further, for each initial conditions 
$({\mathtt x}_0, w_0)$, $(\widetilde{\mathtt x}_0, {\widetilde{w}}_0) \in E \times V$,
there is a constant $c > 0$ such that 
\begin{equation}\label{EST333}
\| {\mathtt x} - {\widetilde{{\mathtt x}}} \|_{H^1(0, T; E)} 
+ \| w - {\widetilde{w}} \|_{L^2(0, T; V)} \le 
c \, (\| {\mathtt x}_0- \widetilde{\mathtt x}_0 \|_E +
\| w_0 - \widetilde {w}_0 \|),
\end{equation}
where $({\mathtt x}, w)$ and $({\widetilde{{\mathtt x}}}, 
{\widetilde{w}})$ 
denote the unique solutions to Problem~$\ref{DVHI}$ 
corresponding to
$({\mathtt x}_0, w_0)$ and $(\widetilde{\mathtt x}_0, {\widetilde{w}}_0)$, respectively. 
Further, the solution has the regularity 
$({\mathtt x}, w) \in C(0,T; E \times H)$. 
\end{Theorem}
\noindent 
{\bf Proof.} \ 
The proof is carried out in several steps and it is 
inspired by a recent result from~\cite[Theorem~20]{MigorskiBiao} combined with a fixed-point argument of Lemma~\ref{CONTR}.

\smallskip

\noindent 
{\bf Step 1}.
Let $\lambda \in L^2(0, T; E)$, 
$\xi \in L^2(0,T; V^*)$, 
$\zeta \in L^2(0, T; Z)$, and
$\eta \in L^2(0, T; Y)$ be fixed
and consider the following auxiliary problem.
\begin{Problem}\label{ProblemAUX1}
Find $w = w_{\lambda\xi\zeta\eta} \in {\mathbb{W}}$
with $w(t) \in K$ for a.e. $t \in (0,T)$ such that
\begin{equation}\label{Problem1b}
\begin{cases}
\displaystyle
\langle w'(t) + A(t, \lambda(t), w(t)) + \xi(t), 
v - w(t) \rangle
+ \, j^0 (t, \lambda(t),\zeta(t), Mw(t); Mv - Mw(t)) \\[1mm]
\ \ \quad \qquad + \, \varphi (t, \lambda (t),\eta(t), Mv) - \varphi(t, \lambda(t), \eta(t), Mw(t)) \ge 
\langle f(t, \lambda(t)), v-w(t) \rangle
\\[1mm]
\qquad\qquad\qquad \qquad
\ \ \mbox{\rm for all} \ v \in K, 
\ \mbox{\rm a.e.} \ t \in (0,T), \\
w(0) = w_0 .
\end{cases}
\end{equation}
\end{Problem}
\noindent 
We define operator 
${\widetilde{A}} \colon (0,T) \times V \to V^*$, 
functions 
${\widetilde{j}} \colon (0, T) \times X \to \real$, 
${\widetilde{\varphi}} \colon (0, T) \times X 
\to \real$, and 
${\widetilde{f}} \colon (0, T) \to V^*$
by 
\begin{eqnarray*}
&&
{\widetilde{A}} (t, v) = A(t, \lambda(t), v)
\ \ \mbox{for} \ \ v \in V, \, \mbox{a.e.} 
\ \, t \in (0, T), \\[1mm]
&& \quad 
{\widetilde{j}}(t,v) 
= j(t, \lambda(t), \zeta(t), v) 
\ \ \mbox{for} \ \ v \in X, \, \mbox{a.e.} 
\ \, t \in (0, T), \\[1mm]
&&\qquad 
{\widetilde{\varphi}} (t, v) 
= \varphi(t, \lambda(t), \eta(t), v) 
\ \ \mbox{for} \ \ v \in X, \, \mbox{a.e.} 
\ \, t \in (0, T), \\[1mm]
&&\qquad\quad
{\widetilde{f}}(t) = f(t, \lambda(t)) - \xi(t)
\ \ \mbox{for a.e.} \ \, t \in (0, T).
\end{eqnarray*}
With this notation 
Problem~\ref{ProblemAUX1} can be equivalently 
reformulated as follows. 
Find element 
$w = w_{\lambda\xi\zeta\eta} \in {\mathbb{W}}$
with $w(t) \in K$ for a.e. $t \in (0,T)$ such that
\begin{equation}\label{ProblemEQUIV}
\begin{cases}
\displaystyle
\langle w'(t) + {\widetilde{A}}(t, w(t)), v - w(t) \rangle 
+ \, {\widetilde{j}}^0 (t, Mw(t); Mv - Mw(t)) 
\\[1mm]
\quad \qquad 
+ \, {\widetilde{\varphi}} (t, Mv) 
- {\widetilde{\varphi}} (t, Mw(t)) \ge 
\langle {\widetilde{f}}(t), v - w(t) \rangle 
\ \ \mbox{\rm for all} \ v \in K, 
\ \mbox{\rm a.e.} \ t \in (0,T), 
\\[1mm]
w(0) = w_0 .
\end{cases}
\end{equation}

\noindent
We shall establish the properties of the data in 
Problem~\ref{ProblemEQUIV}.

By hypothesis $H(A)$(a), (b), it is clear that 
$A(\cdot,\cdot, v)$ is a Carath\'eodory function 
for all $v\in V$. Therefore, the measurability of $t\mapsto\lambda(t)$ entails, 
see~\cite[Corollary~2.5.24]{DMP2}, that 
${\widetilde{A}}(\cdot, v)$ is measurable 
for all $v\in V$.
We use $H(A)$(c), (d) to get that 
${\widetilde{A}}(t, \cdot)$ is demicontinuous for a.e. $t \in (0, T)$ and 
$$
\| {\widetilde{A}}(t, v) \|_{V^*}
\le {\widetilde{a}}_0(t) + a_2 \| v \|
\ \ \mbox{for all} \ \ v \in V, 
\ \mbox{a.e.} \ \, t \in (0, T), 
$$
where ${\widetilde{a}}_{0} \in L^2(0, T)$. 
Further, we employ $H(A)$(d) to infer that
${\widetilde{A}}(t, \cdot)$ is strongly monotone with constant $m_A > 0$ for a.e. $t \in (0, T)$.
The function ${\widetilde{j}}(\cdot, v)$
is measurable because of $H(j)$(c) and of 
the measurability 
of $t\mapsto \lambda(t)$ and 
$t\mapsto \zeta(t)$. It is clear that 
${\widetilde{j}}(t,\cdot)$ is locally Lipschitz for a.e. $t \in (0,T)$ and
$$
\| \partial {\widetilde{j}} (t, v)\|_{X^*} 
\le {\widetilde{c}}_{0j} + 
c_{3j} \| v \|_X
\ \ \mbox{for all} \ \ v \in X, 
\ \mbox{a.e.} \ \, t \in (0, T).
$$
with ${\widetilde{c}}_{0j} \in L^2(0, T)$.
From $H(j)$(e), we know that 
$$
{\widetilde{j}}^0(t, v_1; v_2-v_1) + 
{\widetilde{j}}^0 (t, v_2; v_1-v_2) \le 
m_j \, \| v_1 - v_2 \|_X^2
\ \ \mbox{for all} \ \ v_1, v_2 \in X, 
\ \mbox{a.e.} \ \, t \in (0, T).
$$
By $H(\varphi)$, we deduce that 
${\widetilde{\varphi}}(\cdot, v)$ is measurable 
for all $v \in X$, 
${\widetilde{\varphi}}(t, \cdot)$ is convex and lower 
semicontinuous for a.e. $t \in (0,T)$, and
$$
\| \partial {\widetilde{\varphi}}(t, v) \|_{X^*}
\le 
{\widetilde{c}}_{0\varphi}(t) + c_{3\varphi} \| v \|_X
\ \ \mbox{for all} \ \ v \in X, 
\ \mbox{a.e.} \ \, t \in (0, T)
$$
with ${\widetilde{c}}_{0\varphi} \in L^2(0,T)$.
Since $\lambda \in L^2(0, T; E)$ and 
$\xi \in L^2(0,T; V^*)$, we use $H(f)$ to obtain 
${\widetilde{f}} \in L^2(0,T; V^*)$.
Having verified the above properties of the data, 
and having in mind hypotheses 
$H(K)$, $(H_1)$--$(H_3)$, 
we are in a position to apply~\cite[Theorem~20]{MigorskiBiao} to deduce that 
the problem (\ref{ProblemEQUIV}), and equivalently Problem~\ref{ProblemAUX1}, is uniquely solvable.

\medskip

\noindent
{\bf Step 2}. 
Let $(\lambda_i, \xi_i, \zeta_i, \eta_i) \in
L^2(0,T;E \times V^* \times Z \times Y)$, $i=1$, $2$
and
$w_1 = w_{\lambda_1\xi_1\zeta_1\eta_1}$,
$w_2 = w_{\lambda_2\xi_2\zeta_2\eta_2} \in {\mathbb{W}}$
with $w_1(t)$, $w_2(t) \in K$
for a.e. $t \in (0,T)$,
be the unique solutions to (\ref{ProblemEQUIV}) 
corresponding to 
$(\lambda_1,\xi_1, \zeta_1, \eta_1)$ and
$(\lambda_2, \xi_2, \zeta_2, \eta_2)$, respectively.
We shall prove the following estimate
\begin{eqnarray}\label{888}
&&
\| w_1 - w_2 \|_{L^2(0, t; V)} \le c \,
\big( 
\| \lambda_1-\lambda_2\|_{L^2(0,t;E)} 
+ \| \xi_1 - \xi_2 \|_{L^2(0, t; V^*)}
\nonumber
\\[1mm]
&&
\qquad \qquad 
+ \| \zeta_1 - \zeta_2 \|_{L^2(0, t; Z)} +
\| \eta_1 - \eta_2 \|_{L^2(0, t; Y)}  \big)
\end{eqnarray}
for all $t \in [0,T]$, where $c >0$ is a constant.
From Problem~\ref{ProblemAUX1} it follows that
\begin{eqnarray*}
&&\hspace{-0.6cm}
\displaystyle
\langle w_1'(t) + A(t, \lambda_1(t), w_1(t)) 
+ \xi_1 (t), w_2(t) - w_1(t) \rangle
\\[1mm]
&&\hspace{-0.6cm}
\quad 
+ \, j^0 (t, \lambda_1(t),\zeta_1(t), Mw_1(t); 
Mw_2(t) - Mw_1(t)) \\[1mm]
&&\hspace{-0.6cm}
\qquad 
+ \, \varphi (t, \lambda_1 (t), \eta_1(t), Mw_2(t)) - \varphi(t, \lambda_1(t), \eta_1(t), Mw_1(t)) \ge 
\langle f(t, \lambda_1(t)), w_2(t)-w_1(t) \rangle	
\end{eqnarray*}
for a.e. $t \in (0, T)$ 
and
\begin{eqnarray*}
	&&\hspace{-0.6cm}
	\displaystyle
	\langle w_2'(t) + A(t, \lambda_2(t), w_2(t)) 
	+ \xi_2 (t), w_1(t) - w_2(t) \rangle
	\\[1mm]
	&&\hspace{-0.6cm}
	\quad 
	+ \, j^0 (t, \lambda_2(t),\zeta_2(t), Mw_2(t); 
	Mw_1(t) - Mw_2(t)) \\[1mm]
	&&\hspace{-0.6cm}
	\qquad 
	+ \, \varphi (t, \lambda_2 (t), \eta_2(t), Mw_1(t)) - \varphi(t, \lambda_2(t), \eta_2(t), Mw_2(t)) \ge 
	\langle f(t, \lambda_2(t)), w_1(t)-w_2(t) \rangle	
\end{eqnarray*}
for a.e.\ $t \in (0, T)$, and
$w_1(0) - w_2 (0) = 0$.
Next, we sum up the last two inequalities to get
\begin{eqnarray*}
&&\hspace{-0.7cm}
\langle w'_1(t) - w'_2(t), w_1(t) - w_2(t) 
\rangle +
\langle A(t, \lambda_1(t), w_1(t)) 
- A(t, \lambda_2(t), w_2(t)), w_1(t) - w_2(t)\rangle 
\\ [1mm]
&&\hspace{-0.5cm}
\le 
\langle \xi_1(t) -\xi_2(t), w_2(t) - w_1(t) \rangle
+ 
\langle f(t, \lambda_1(t)) - f(t, \lambda_2(t)), 
w_1(t) - w_2(t) \rangle
\\[1mm]
&&\hspace{-0.3cm}
+ j^0(t, \lambda_1(t), \zeta_1(t), Mw_1(t); Mw_2(t) - Mw_1(t)) 
+ j^0(t, \lambda_2(t), \zeta_2(t), Mw_2(t); Mw_1(t) - Mw_2(t)) \\ [1mm]
&&\hspace{-0.1cm}
+ \, \varphi(t, \lambda_1(t), \eta_1(t), Mw_2(t)) 
- \varphi (t, \lambda_1(t), \eta_1(t), Mw_1 (t))
\\[1mm]
&&\hspace{0.1cm}
+ \, \varphi(t, \lambda_2(t), \eta_2(t), Mw_1(t)) 
- \varphi (t, \lambda_2(t), \eta_2(t), Mw_2 (t)) 
\end{eqnarray*}
for a.e. $t \in (0,T)$. 
We integrate the above inequality on $(0, t)$, 
apply the integration by parts formula, 
see~\cite[Proposition~3.4.14]{DMP2}, 
to the first term, 
use hypotheses $H(A)$(d), $H(j)$(e), $H(\varphi)$(e) 
to obtain
\begin{eqnarray*}
&&\hspace{-0.5cm}
\frac{1}{2}\| w_1(t) - w_2(t) \|^2_H 
- \frac{1}{2}\| w_1(0) - w_2(0) \|^2_H +
m_A \, \int_{0}^{t} \| w_1(s) - w_2(s) \|^2 \, ds 
\\
&&\hspace{-0.3cm}
\le 
{\bar{m}}_A \int_0^t
\| \lambda_1(s)-\lambda_2(s)\|_E \, 
\| w_1(s)-w_2(s)\| \, ds
+ \int_0^t 
\| \xi_1(s)-\xi_2(s)\|_{V^*} \, \|w_1(s)-w_2(s)\| \, ds
\\
&&\hspace{-0.1cm}
+ \int_0^t 
\| f(s, \lambda_1(s))-f(s, \lambda_2(s)) \|_{V^*}
 \, \| w_1(s) - w_2(s) \| \, ds
+ m_j \int_0^t \| Mw_2(s) - Mw_1(s) \|_X^2
\\
&&\hspace{0.1cm}
+ \, {\bar{m}}_{j} \int_{0}^{t} 
\Big(
\| \lambda_1(s)-\lambda_2(s)\|_E 
+ \| \zeta_1 (s) - \zeta_2(s) \|_Z 
\Big) 
\| w_1(s) - w_2(s)\| \, ds \\
&&\hspace{0.3cm}
+ \, m_\varphi \int_{0}^{t} 
\Big(
\| \lambda_1(s)-\lambda_2(s)\|_E 
+ \| \eta_1 (s) - \eta_2(s) \|_Y 
\Big) \| w_1(s) - w_2(s)\| \, ds.
\end{eqnarray*}
for all $t \in [0,T]$. 
Next, using hypothesis $H(f)$(b), 
condition $w_1(0) - w_2(0) = 0$ and 
the inequality 
$\| M v \|_X \le \| M_1\|\, \| v \|$
for all $v \in V$, 
and the H\"older inequality, we have
\begin{eqnarray*}
&&
\big( m_A - m_j \|M_1\|^2 \big) \,
\| w_1 - w_2 \|_{L^2(0, t; V)}^2
\le 
\| \xi_1 - \xi_2 \|_{L^2(0, t; V^*)}
\| w_1 - w_2 \|_{L^2(0, t; V)}
\\[2mm]
&&\quad 
+ \, 
\Big(
{\bar{m}}_A + L_f+{\bar{m}}_j \| M_1\| 
+ m_\varphi \| M_1\|
\Big) \,
\| \lambda_1 - \lambda_2 \|_{L^2(0, t; E)}
\| w_1 - w_2 \|_{L^2(0, t; V)} 
\\[2mm] 
&&\qquad 
+ \, 
\Big( 
{\bar{m}}_j \| M_1\| \,
\| \zeta_1 - \zeta_2 \|_{L^2(0, t; Z)}
+ m_\varphi \| M_1\| 
\| \eta_1 - \eta_2 \|_{L^2(0, t; Y)} 
\Big)
\, \| w_1 - w_2 \|_{L^2(0, t; V)}
\end{eqnarray*}
for all $t \in [0,T]$.
Hence, by $(H_1)$, the inequality (\ref{888}) follows.

\medskip

\noindent 
{\bf Step 3.}
Let $w \in L^2(0, T; V)$ be fixed. 
We claim that there exists a unique function 
${\mathtt x} \in H^1(0, T; E)$ 
solution to the Cauchy problem
\begin{eqnarray}\label{ODE}
\begin{cases}
\displaystyle
{\mathtt x}'(t) = F(t, {\mathtt x}(t), w(t), (Sw)(t)) 
\ \ \mbox{for a.e.} \ \ t \in (0, T)
\\
{\mathtt x}(0) = {\mathtt x}_0.
\end{cases}\end{eqnarray}
Further, if 
$w_1$, $w_2 \in L^2(0, T; E)$
and ${\mathtt x}_1$, ${\mathtt x}_2 \in H^1(0,T;E)$ denote the unique solutions to (\ref{ODE}) 
corresponding to $w_1$ and $w_2$, respectively, then
\begin{equation}\label{EST1}
\| {\mathtt x}_1(t) - {\mathtt x}_2(t)\|_E \le 
m \int_0^t \| w_1(s) - w_2(s)\| \, ds
\ \ \mbox{for all} \ \ t \in (0, T)
\end{equation}
with a constant $m > 0$.
In fact, let $w \in L^2(0, T; V)$ 
and define the function
${\widetilde{F}} \colon (0, T) \times E \to E$ by
$$
{\widetilde{F}} (t, {\mathtt x}) 
= F(t, {\mathtt x}, w(t), (S w)(t)) 
\ \ \mbox{for} \ \ {\mathtt x} \in E, 
\ \mbox{a.e.} \ \ t \in (0, T).
$$
By hypotheses $H(F)$, we obtain
$$
\| {\widetilde{F}}(t, {\mathtt x}_1) - {\widetilde{F}}(t, {\mathtt x}_2)\|_E
= 
\| F(t, {\mathtt x}_1, w(t), (S w)(t)) 
- F(t, {\mathtt x}_2, w(t), (S w)(t)) \|_E 
\le L_f \| {\mathtt x}_1 - {\mathtt x}_2\|_E
$$
for all ${\mathtt x}_1$, ${\mathtt x}_2 \in E$, 
a.e. $t\in (0, T)$, 
and
$$
(0, T) \ni t \mapsto {\widetilde{F}}(t, {\mathtt x}) 
\in E
$$
is a $L^2$-function for any ${\mathtt x} \in E$.
Hence, we are able to apply the classical Cauchy-Lipschitz theorem, see e.g.,~\cite[Theorem~2.30]{SST} to deduce the unique solvability of (\ref{ODE}). 
Next, for $w_1$, $w_2 \in L^2(0,T; V)$, 
we integrate equation in (\ref{ODE}) to get 
$$
{\mathtt x}_i(t) = {\mathtt x}_0 
+ \int_0^t 
F(s, {\mathtt x}_i(s), w_i(s), (S w_i)(s))\, ds 
\ \ \mbox{for} \ \ i=1,2.
$$
Thus
\begin{eqnarray*}
&&\hspace{-0.6cm}
\| {\mathtt x}_1(t) -{\mathtt x}_2(t) \|_E
\le 
\int_0^t 
\|F(s, {\mathtt x}_1(s), w_1(s), (S w_1)(s))
-F(s, {\mathtt x}_2(s), w_2(s), (S w_2)(s))
\|_E \, ds \\[1mm]
&&\hspace{-0.6cm}
\quad \le 
L_F \int_0^t \|{\mathtt x}_1(s)-{\mathtt x}_2(s)\|_E 
+ \| w_1(s)-w_2(s)\| + 
\| (S w_1)(s)-(S w_2)(s)\|_U 
\, ds
\\[1mm]
&&
\hspace{-0.6cm}
\qquad \le 
L_F \int_0^t \|{\mathtt x}_1(s)-{\mathtt x}_2(s)\|_E 
+ \| w_1(s)-w_2(s)\| \, ds + 
c_S L_F T \int_0^t \| w_1(s)-w_2(s)\| 
\, ds
\end{eqnarray*}
for all $t \in [0, T]$. 
Finally, we use the standard argument based on the Growall inequality to deduce inequality (\ref{EST1}) 
with a constant $m>0$.  

We denote by 
$R_0 \colon L^2(0,T;V) \to H^1(0,T;E) \subset L^2(0,T;E)$
the operator $R_0 (w) = {\mathtt x}$ which to any 
$w \in L^2(0,T;V)$ assigns the unique solution 
${\mathtt x} \in L^2(0,T;E)$ 
to the problem (\ref{ODE}). 
The inequality (\ref{EST1}) can be restated as 
follows 
\begin{equation}\label{EST2}
\| (R_0 w_1)(t) - (R_0 w_2)(t)\|_E \le 
m \int_0^t \| w_1(s) - w_2(s)\| \, ds
\ \ \mbox{for all} \ \ t \in (0, T), 
\end{equation}
i.e., $R_0$ is a history-dependent operator.

\medskip

\noindent
{\bf Step 4}.
In this part of the proof we apply a fixed point argument. We define the operator
$\Lambda \colon L^2(0,T; E \times V^* \times Z \times Y) \to L^2(0,T;E \times V^* \times Z \times Y)$ by
$$
\Lambda(\lambda, \xi, \zeta, \eta) =
(R_0 w_{\lambda\xi\zeta\eta},
{R}_1 w_{\lambda\xi\zeta\eta}, 
{R_2} w_{\lambda\xi\zeta\eta}, 
{R_3} w_{\lambda\xi\zeta\eta})
$$
for all 
$(\lambda,\xi, \zeta, \eta) \in L^2(0,T;E \times V^* \times Z \times Y)$, 
where $w_{\lambda\xi\zeta\eta} \in {\mathbb{W}}$ 
denotes the unique solution
to Problem~\ref{ProblemAUX1} corresponding to 
$(\lambda, \xi, \zeta, \eta)$.
From hypothesis $(H_3)$(b)--(d), inequalities (\ref{888}) and (\ref{EST2}), and by the H\"older inequality, 
we find a constant $c > 0$ such that
\begin{eqnarray*}
&&\hspace{-0.6cm}
\| \Lambda(\lambda_1,\xi_1,\zeta_1,\eta_1)(t) -
\Lambda(\lambda_1, \xi_2, \zeta_2,\eta_2)(t) 
\|^2_{E\times V^* \times Z \times Y} \\[2mm]
&&
\hspace{-0.4cm} =
\| (R_0 w_1)(t) - (R_0 w_2)(t) \|_{E}^2
+ \| ({R}_1 w_1)(t) - ({R}_1 w_2)(t) \|_{V^*}^2
\\[2mm]
&&\hspace{-0.2cm}
+ \| ({R_2} w_1)(t) - ({R_2} w_2)(t) \|_Z^2
+ \| ({R_3} w_1)(t) - ({R_3} w_2)(t) \|_Y^2 \\[1mm]
&&
\hspace{-0.0cm} \le 
\Big( m \int_0^t \| w_1(s) - w_2(s) \| \, ds \Big)^2
+ \Big( c_{R_1} \int_0^t \| w_1(s) - w_2(s) \| \, ds \Big)^2
\\[2mm]
&&
\hspace{0.2cm}
+ \Big( c_{R_2} \int_0^t \| w_1(s) - w_2(s) \| \, ds \Big)^2 
+ \Big( c_{R_3} \int_0^t \| w_1(s) - w_2(s) \| \, ds \Big)^2
\le c \, \| w_1 - w_2 \|_{L^2(0, t; V)}^2 \\
&&\hspace{0.4cm}
\le c \, \Big(	
\| \lambda_1-\lambda_2\|^2_{L^2(0,t;E)} 
+ \| \xi_1 - \xi_2 \|^2_{L^2(0, t; V^*)} 
+ \| \zeta_1 - \zeta_2 \|^2_{L^2(0, t; Z)} 
+ \| \eta_1 - \eta_2 \|^2_{L^2(0, t; Y)} \Big),
\end{eqnarray*}
which implies
\begin{eqnarray}
&&\| \Lambda(\lambda_1,\xi_1, \zeta_1, \eta_1)(t) -
\Lambda(\lambda_2, \xi_2, \zeta_2, \eta_2)(t)
\|^2_{E\times V^* \times Z \times Y}
\nonumber \\
&&\qquad\qquad
\le c \, \int_{0}^{t} \, 
\| (\lambda_1, \xi_1, \zeta_1, \eta_1)(s) -
(\lambda_2, \xi_2, \zeta_2, \eta_2)(s) 
\|^2_{E \times V^* \times Z \times Y} \, ds
\label{pstaly}
\end{eqnarray}
for a.e.\ $t \in (0,T)$.
By Lemma~\ref{CONTR},
we deduce that there exists a unique fixed point $(\lambda^*, \xi^*,\zeta^*,\eta^*)$ of $\Lambda$, 
that is 
$$
(\lambda^*, \xi^*, \zeta^*,\eta^*) \in 
L^2(0,T;E \times V^* \times Z \times Y)
\ \ {\rm and}\ \ 
\Lambda(\lambda^*, \xi^*, \zeta^*, \eta^*) 
= (\lambda^*, \xi^*, \zeta^*, \eta^*).
$$
%
%

Let $(\lambda^*, \xi^*, \zeta^*, \eta^*) \in 
L^2(0,T;E \times V^* \times Z \times Y)$ 
be the unique fixed
point of the operator $\Lambda$.
We define $w_{\lambda^* \xi^* \zeta^* \eta^*} 
\in {\mathbb{W}}$ 
to be the unique solution to Problem~\ref{ProblemAUX1} corresponding to $(\lambda^*, \xi^*,\zeta^*, \eta^*)$.
From the definition of the operator $\Lambda$, we have
$$
\lambda^* = R_0(w_{\lambda^* \xi^* \zeta^* \eta^*}), 
\ \ 
\xi^* = {R}_1(w_{\lambda^* \xi^* \zeta^* \eta^*}), 
\ \
\eta^* = {R_2}(w_{\lambda^* \xi^* \zeta^* \eta^*})
\ \ \mbox{and} \ \
\zeta^* = {R_3}(w_{\lambda^* \xi^* \zeta^* \eta^*}) .
$$
Finally, we use these relations in 
Problem~\ref{ProblemAUX1}, and 
conclude that 
$w_{\lambda^* \xi^* \zeta^* \eta^*}\in {\mathbb{W}}$
is the unique solution to Problem~\ref{ProblemAUX1}.

\smallskip

\noindent
{\bf Step 5.}
Finally, we shall prove (\ref{EST333}). Let 
$({\mathtt x}, w)$ and $({\widetilde{{\mathtt x}}}, {\widetilde{w}})$ 
be the unique solutions to Problem~$\ref{DVHI}$ 
corresponding to
$({\mathtt x}_0, w_0)$, $(\widetilde{\mathtt x}_0, {\widetilde{w}}_0) \in E \times V$, respectively.
We apply the arguments used in the inequalities in 
Steps~2 and~3.
We write two inequalities for $w$ and ${\widetilde{w}}$ 
as in (\ref{ProblemEQUIV}), then we choose $v = {\widetilde{w}}(t)$ 
and $v = w(t)$ in the inequality satisfied by 
$w$ and ${\widetilde{w}}$, respectively, 
and add the resulting inequalities to get
\begin{eqnarray*}
	&&\hspace{-0.7cm}
	\langle w'(t) - {\widetilde{w}}'(t), w(t) - {\widetilde{w}}(t) 
	\rangle +
	\langle A(t, \lambda_1(t), w(t)) 
	- A(t, \lambda_2(t), {\widetilde{w}}(t)), 
	w(t) - {\widetilde{w}}(t)\rangle 
	\\ [1mm]
	&&\hspace{-0.5cm}
	\le 
	\langle (R_1w)(t) - (R_1{\widetilde{w}})(t), 
	{\widetilde{w}}(t) - w(t) \rangle +	
	\langle f(t,{\mathtt x}(t)) 
	- f(t, {\widetilde{{\mathtt x}}}(t)),
	{\widetilde{w}}(t) - w(t) \rangle
	\\[1mm]
	&&\hspace{-0.3cm}
	+ \, j^0(t, {\mathtt x}(t), (R_2w)(t); M{\widetilde{w}}(t) 
	- M w(t)) 
	+ j^0(t, {\widetilde{{\mathtt x}}}(t), (R_2{\widetilde{w}})(t); Mw(t) - M{\widetilde{w}}(t)) \\ [1mm]
	&&\hspace{-0.1cm}
	+ \, \varphi(t, {\mathtt x}(t), (R_3w)(t), M{\widetilde{w}}(t)) 
	- \varphi (t, {\mathtt x}(t), (R_3w)(t), Mw (t))
	\\[1mm]
	&&\hspace{0.1cm}
	+ \, \varphi(t, {\widetilde{{\mathtt x}}}(t), 
	(R_3{\widetilde{w}})(t), Mw(t)) 
	- \varphi (t, {\widetilde{{\mathtt x}}}(t), 
	(R_3{\widetilde{w}})(t), M{\widetilde{w}}(t)) 
\end{eqnarray*}
for a.e. $t \in (0,T)$. 
Integrating the latter on $(0, t)$, 
we use the integration by parts formula 
to the first term, 
hypotheses $H(A)$(d), $H(j)$(e), $H(\varphi)$(e), 
$H(f)$(b), and the Young inequality 
$ab \le \frac{\varepsilon^2}{2} a^2 +
\frac{1}{2\varepsilon^2}b^2$ for $a$, $b \in \real$
to obtain
\begin{eqnarray*}
	&&\hspace{-0.5cm}
	(m_A - m_j \| M_1 \|^2)\, 
	\int_{0}^{t} \| w(s) - {\widetilde{w}}(s) \|^2 \, ds
	\le \frac{1}{2}\| w_0 - \widetilde w_0 \|^2_H 
	\\
	&&\hspace{-0.3cm}
	+ \,  
	\frac{\varepsilon^2}{2} \left(
	{\bar{m}}_A + 2 {\bar{m}}_j \| M_1\| + 2m_\varphi 
	\| M_1\| +2 \right)
	\int_{0}^{t} \| w(s) - {\widetilde{w}}(s) \|^2 \, ds \\
	&&
	+ \, \frac{1}{2\varepsilon^2}
	\left(	{\bar{m}}_A + {\bar{m}}_j \| M_1 \| 
	+ m_\varphi \| M_1 \| + L_f^2 \right)
	\int_{0}^{t} \| {\mathtt x}(s)
	-{\widetilde{{\mathtt x}}}(s) \|^2 \, ds
	\\
	&&
	+ \, \frac{1}{2\varepsilon^2} 
	\int_{0}^{t} \| (R_1w)(s)-(R_1{\widetilde{w}})(s) \|_{V^*}^2 \, ds
	\\
	&&
	+ \, \frac{1}{2\varepsilon^2} {\bar{m}}_j \| M_1\| 
	\int_{0}^{t} \| (R_2w)(s)-(R_2{\widetilde{w}})(s) \|_{Z}^2 \, ds
	\\
	&&
	+ \, \frac{1}{2\varepsilon^2} m_\varphi \| M_1\|
	\int_{0}^{t} \| (R_3w)(s)-(R_3{\widetilde{w}})(s) \|_{V^*}^2 \, ds.
\end{eqnarray*}
for all $t \in [0,T]$. 
By the smallness condition $(H_1)$, we choose 
$\varepsilon > 0$ such that 
$$
m_A - m_j \| M_1 \|^2 - \frac{\varepsilon^2}{2} \left(
{\bar{m}}_A + 2 {\bar{m}}_j \| M_1\| + 2m_\varphi 
\| M_1\| +2 \right) > 0.
$$
Next, we employ $(H_3)$(b)--(d) and get
\begin{eqnarray}
&&\label{EST88}
d_1 \int_{0}^{t} \| w(s) - {\widetilde{w}}(s) \|^2 \, ds 
\le \frac{1}{2}\| w_0 - \widetilde w_0 \|^2_H 
\\
&&\quad 
+ \, d_2
\int_{0}^{t} \| {\mathtt x}(s)-{\widetilde{{\mathtt x}}}(s) \|_E^2 
\, ds
+ \, d_3 \int_0^t 
\left( \int_{0}^{s} \| w(\tau) - {\widetilde{w}}(\tau) \|^2 \, d\tau
\right)\, ds
\nonumber 
\end{eqnarray}
with some positive constants $d_i$, $i=1$, $2$, $3$.
On the other hand, we use the ordinary differential 
equation in Problem~\ref{DVHI} and, analogously as in (\ref{EST1}), by the Gronwall inequality, we have
\begin{equation}\label{EST99}
\| {\mathtt x}(t) - {\widetilde{{\mathtt x}}}(t)\|_E 
\le 
d_4 
\left(
\| {\mathtt x}_0 - \widetilde{\mathtt x}_0\|_E
+ \int_0^t \| w(s) - {\widetilde{w}}(s)\| \, ds
\right)
\ \ \mbox{for all} \ \ t \in (0, T)
\end{equation}
with a constant $d_4 > 0$.
We combine (\ref{EST88}) and (\ref{EST99}), and again use the Gronwall inequality to deduce 
\begin{equation}\label{EST901}
\| w - {\widetilde{w}} \|_{L^2(0, t; V)} \le 
d_5 \left( 
\| {\mathtt x}_0 - \widetilde{\mathtt x}_0\|_E 
+ \| w_0 - \widetilde w_0 \| \right) 
\end{equation}
for all $t \in [0, T]$ with $d_5 > 0$.
By $H(F)$(b), it follows
\begin{equation*}\label{INEQ902}
\| {\mathtt x}' - {\widetilde{{\mathtt x}}}'\|_{L^2(0, t; E)} 
\le d_6 \left( 
\| {\mathtt x} - {\widetilde{{\mathtt x}}}\|_{L^2(0, t; E)}
+ \| w - {\widetilde{w}}\|_{L^2(0, t; V)} \right)
\ \ \mbox{for all} \ \ t \in [0, T]
\end{equation*}
with a constant $d_6 > 0$.
The latter together with (\ref{EST99}) and (\ref{EST901}) imply (\ref{EST333}). 
The regularity of the solution follows from the
following embeddings 
$W^{1,2}(0, T; E) \subset C([0, T]; E)$ and 
${\mathbb{W}} \subset C([0,T]; H)$.
This completes the proof of the theorem.
\hfill$\Box$

We conclude this section with a corollary of  Theorem~\ref{Theorem1} when the term $f$ 
is independent of the second variable 
and $f \in L^2(0, T; V^*)$.
\begin{Corollary}\label{Corollary2}
Let hypotheses $H(F)$, $H(A)$, $H(j)$, $H(\varphi)$, $H(K)$, $H(M)$, $(H_1)$--$(H_3)$ hold, 
and $f \in L^2(0, T; V^*)$. 
Then Problem~$\ref{DVHI}$ is uniquely solvable with  
$({\mathtt x}, w) \in H^1(0,T; E) \times {\mathbb{W}}$ 
and $w(t) \in K$ for a.e. $t \in (0, T)$.
Moreover, for every 
$({\mathtt x}_0, w_0, f)$, 
$(\widetilde{\mathtt x}_0, {\widetilde{w}}_0, {\widetilde{f}}) \in E \times V \times L^2(0, T; V^*)$,
there exists a constant $c > 0$ such that 
	\begin{equation}\label{EST444}
	\| {\mathtt x} - {\widetilde{{\mathtt x}}} \|_{H^1(0, T; E)} 
	+ \| w - {\widetilde{w}} \|_{L^2(0, T; V)} \le 
	c \, (\| {\mathtt x}_0- \widetilde{\mathtt x}_0 \|_E 
	+ \| w_0 - \widetilde {w}_0\|
	+ \| f - {\widetilde{f}} \|_{L^2(0, T; V^*)}),
	\end{equation}
	where $({\mathtt x}, w)$ and 
	$({\widetilde{{\mathtt x}}}, {\widetilde{w}})$ 
	are the unique solutions to Problem~$\ref{DVHI}$ 
	corresponding to
	$({\mathtt x}_0, w_0, f)$ and 
	$(\widetilde{\mathtt x}_0, {\widetilde{w}}_0, \widetilde{f})$, respectively.
\end{Corollary}

\section{A dynamic frictionless viscoplastic contact problem}\label{Application1}

In this section we shall illustrate the applicability 
of results of Section~\ref{s1}. We examine an example of 
a dynamic unilateral viscoplastic frictionless 
contact problem.
We describe the physical setting and give the classical 
formulation of the contact problem, provide its variational formulation, and prove a result on its well-posedness.

A viscoplastic body occupies a bounded domain $\Omega$ of $\mathbb{R}^{d}$, $d=1$, $2$, $3$.  
The boundary $\Gamma$ of $\Omega$ is supposed to be Lipschitz continuous and splitted into three mutually disjoint and measurable parts 
$\Gamma= \overline{\Gamma}_{D} \cup 
\overline{\Gamma}_{N} \cup 
\overline{\Gamma}_{C}$ with
$|\Gamma_{D}| > 0$. The outward unit normal 
at $\Gamma$, denoted by $\bnu$, is defined a.e. 
on $\Gamma$.
Let $Q = \Omega \times (0, T)$, where $0 < T < \infty$.

The symbol $\mathbb{S}^{d}$ denotes the
space of $d\times d$ symmetric matrices. 
We use the standard notation for inner products 
and norms on $\mathbb{R}^{d}$ and $\mathbb{S}^{d}$, i.e., 
$\bu\cdot \bv=u_{i}v_{i}$, 
$\|\bu\|^2=\bu\cdot \bu$ for $\bu=(u_{i})$,
$\bv=(v_{i})\in \mathbb{R}^{d}$, and
$\bsigma\cdot \btau=\sigma_{i}\tau_{i}$,
$\|\bsigma\|^2= \bsigma\cdot \bsigma$ for 
$\bsigma= (\sigma_{ij})$, 
$\btau=(\tau_{ij})\in \mathbb{S}^{d}$,
respectively. 
The normal and tangential components of 
vectors and tensors at the boundary 
are expressed by the following notation
$$
v_{\nu}=\bv\cdot \bnu, \ \ \ 
\bv_{\tau}=\bv-v_{\nu}\, \bnu, \ \ \ 
\sigma_{\nu}=(\bsigma\bnu)\cdot\bnu, \ \ \ 
\bsigma_{\tau}=\bsigma\, \bnu - \sigma_{\nu}\, \bnu.
$$
Since $\bv_\tau \cdot \bnu = \bsigma_\tau \cdot \bnu = 0$, 
the following decomposition formula holds
\begin{equation}\label{DEC}
\bsigma\bnu \cdot \bv 
= (\sigma_{\nu} \, \bnu + \bsigma_{\tau})
\cdot 
(v_\nu \, \bnu + \bv_\tau) = 
\sigma_{\nu} \, v_\nu + \bsigma_{\tau} \cdot \bv_{\tau}.
\end{equation}

The classical model for the contact process has the following formulation. 
\begin{Problem}\label{Plastic}
Find a displacement field
$\bu \colon Q \to\mathbb{R}^d$ and a stress field
$\bsigma \colon Q \rightarrow \mathbb{S}^d$ 
such that for all $t\in (0,T)$,
\begin{align}
\label{equation1}
\bsigma(t) &={\mathscr A}\bvarepsilon ({\bu}'(t))
+{\mathscr B}\bvarepsilon({\bu}(t))
+\int_0^t{\mathscr G}\big(
\bsigma(s) - 
{\mathscr A}\bvarepsilon ({\bu}'(s)), 
\bvarepsilon({\bu}(s) \big)\,ds \quad
&{\rm in}\
	&\Omega,\\[1mm]
	\label{equation2}
	{\bu}''(t)&={\rm Div}\,\bsigma(t)+\fb_0(t)\quad&{\rm in}\ &\Omega,\\[1mm]
	\label{equation3} \bu(t)&=\bzero &{\rm on}\ &\Gamma_D,\\[1mm]
	\label{equation4} \bsigma(t)\bnu&=\fb_N(t)\quad&{\rm on}\ &\Gamma_N,\\[1mm]
	\nonumber
	u'_{\nu}(t)&\le g, \
	\sigma_{\nu}(t)+\xi(t) \le 0, \ (u'_{\nu}(t)-g)(\sigma_{\nu}(t)+\xi(t))=0 
	\\[1mm]
	\label{equation5}
	&\ \ \ {\rm with} \ \xi(t) \in k(u_{\nu}(t))\,
	\partial j_{\nu}(u'_{\nu}(t))
	\quad&{\rm	on}\
	&\Gamma_C,\\[1mm]
	\label{equation6} 
	\bsigma_\tau(t) &=\bzero\quad&{\rm on}\ &\Gamma_C,
	%
	\end{align} 
	and
	\begin{equation}\label{equation7}
	\bu(0)=\bu_0,\qquad
	{\bu}'(0)=\bw_0\quad \ {\rm in}\quad \Omega.
	\end{equation}
\end{Problem}

\smallskip

\noindent
In this problem, equation (\ref{equation1}) represents a general viscoplastic constitutive law in which $\mathscr{A}$, $\mathscr{B}$ and $\mathscr{G}$ 
stand for the viscosity operator, the elasticity operator 
and the viscoplastic constitutive function, respectively. 
Recall that the components of the linearized strain tensor 
$\bvarepsilon(\bu)$ are given by
$$
\bvarepsilon_{ij}(\bu) = (\bvarepsilon(\bu))_{ij}
= \frac{1}{2} (u_{i,j} + u_{j,i})
\ \ \mbox{in} \ \ \Omega,
$$
where $u_{i,j} = \partial u_i/\partial x_j$.
Equation (\ref{equation2}) is the equation 
of motion which governs the evolution of the mechanical state of the body. 
Here 
$u''(t) = \partial^2\bu / \partial t^2$
represents the acceleration field, 
${\rm Div} \, \bsigma = (\sigma_{ij,j})$, and
$\fb_{0}$ denotes the density of body forces.
The displacement homogeneous boundary condition
(\ref{equation3}) means that the body is fixed on $\Gamma_{D}$, while (\ref{equation4}) is the traction boundary condition with surface tractions of density $\fb_N$ acting on $\Gamma_{N}$.
Condition (\ref{equation5}) is the Signorini  unilateral contact boundary condition for the normal velocity in which $g > 0$ and $\partial j_\nu$ denotes the Clarke subgradient of a prescribed function $j_\nu$.
Condition $\xi(t) \in k(u_{\nu}(t))\,
\partial j_{\nu}(u'_{\nu}(t))$ on $\Gamma_C$ represents the normal damped response
condition where $k$ is a given damper
coefficient depending on the normal displacement.
Condition (\ref{equation6}) is called the frictionless condition in which the tangential part of the stress vanishes.
%
Examples and details on mechanical interpretation of conditions (\ref{equation5}) and (\ref{equation6}) can be found in~\cite{HMS2017,MOSBOOK}
and references therein.
Finally, equations (\ref{equation7}) are the initial conditions in which $\bu_0$ and $\bw_0$ represent the initial displacement and the initial velocity, respectively.

To cope with the weak formulation of Problem~\ref{Plastic}, 
we introduce the spaces
\begin{equation}\label{SPACESVH}
V=\{\, \bv\in H^{1}(\Omega;\mathbb{R}^{d})\mid \bv={\bf{0}}~\textrm{on}~\Gamma_{D}\, \} 
\ \ \mbox{and} \ \  
\mathcal{H}=L^{2}(\Omega;\mathbb{S}^{d}).
\end{equation}
The inner product and the corresponding norm on $V$
are given by
\begin{equation}\label{NORM}
(\bu,\bv)_{V} = (\bvarepsilon(\bu), \bvarepsilon(\bv))_{\mathcal{H}},~~\|\bv\|_{V}=\| \bvarepsilon(\bv)\|_{\mathcal{H}}
\ \ \mbox{for all} \ \ \bu,\, \bv\in V. 
\end{equation}
The space $\mathcal{H}$ is a Hilbert space endowed with the inner product
\begin{equation*}
\langle \bsigma, \btau \rangle_{\mathcal{H}}=
\int_{\Omega}\bsigma(\bx) \cdot \btau (\bx)\, dx
\ \ \mbox{for} \ \ \bsigma, \btau \in {\mathcal H},
\end{equation*}
and the associated norm $\|\cdot\|_{\mathcal{H}}$.
Since $|\Gamma_D| > 0$, the Korn inequality 
%
implies 
that the norm $\| \cdot \|_V$ defined by (\ref{NORM}) 
is equivalent to the usual 
norm $\| \cdot \|_{H^1 (\Omega;\real^d)}$.
Recall that the trace operator
$\gamma\colon V\to L^{2}(\Gamma_C;\mathbb{R}^{d})$
is linear and continuous, i.e.,
$\|\bv \|_{L^{2}(\Gamma;\mathbb{R}^{d})}
\le \|\gamma\|\|\bv \|_{V}$ for $\bv\in V$, 
where $\|\gamma\|$ denotes the norm of the trace
operator in ${\cal L}(V,L^2(\Gamma_C;\real^d))$.

We need the following hypotheses 
on the data to Problem~\ref{Plastic}.

\smallskip

\noindent
$\underline{H({\mathscr{A}})}:$ \quad
${\mathscr{A}} \colon Q \times \mathbb{S}^{d}\rightarrow \mathbb{S}^{d}$
is such that

\smallskip

\lista{
	\item[(1)]
	${\mathscr{A}}(\cdot,\cdot,\bvarepsilon)$
	is measurable on $Q$ for all
	$\bvarepsilon \in \mathbb{S}^{d}$.
\smallskip
\item[(2)]
${\mathscr{A}}(\bx, t,\cdot)$ 
is continuous on $\mathbb{S}^{d}$ for a.e. $t \in (0, T)$. 
\smallskip	
\item[(3)]
$\| {\mathscr{A}} (\bx, t, \bvarepsilon) \| 
\le {\widetilde{a}}_0(\bx,t) + {\widetilde{a}}_2 \, 
\| \bvarepsilon \|$ for all 
$\bvarepsilon \in \mathbb{S}^{d}$, 
a.e. $(\bx, t) \in Q$ with ${\widetilde{a}}_0$, ${\widetilde{a}}_2 \ge 0$, ${\widetilde{a}}_0 \in L^2(Q)$. 
	\smallskip
	\item[(4)]
	$\mbox{there exists} \ m_{\mathscr{A}}>0 \ \mbox{such that} \ \
	(\mathscr{A}(\bx,t,\bvarepsilon_{1})-
	\mathscr{A}(\bx,t,\bvarepsilon_{2}))\cdot (\bvarepsilon_{1}-\bvarepsilon_{2})
	\ge m_{\mathscr{A}} \, 
	\|\bvarepsilon_{1}-\bvarepsilon_{2}\|^{2}$ \\
	$\ \ \qquad \mbox{for all} \  \bvarepsilon_{1}, \bvarepsilon_{2}\in \mathbb{S}^{d}, \ \mbox{a.e.}
	\ (\bx, t)\in Q$.
	\smallskip
	\item[(5)]
	${\mathscr A}(\bx,t,{\bf{0}}) = {\bf{0}}$
	for a.e. $(\bx,t) \in Q$.
}

\smallskip

\noindent
$\underline{H({\mathscr B})}:$ \quad
$\mathscr {B}\colon Q \times \mathbb{S}^{d}
\rightarrow \mathbb{S}^{d}$
is such that

\smallskip

\lista{
	\item[(1)]
	$\mathscr{B}(\cdot,\cdot,\bvarepsilon)$
	is measurable on $Q$ for all
	$\bvarepsilon \in \mathbb{S}^{d}$.
	\smallskip
	\item[(2)]
	there exists $L_{\mathscr{B}}>0$ such that
	$\|\mathscr{B}(\bx,t, \bvarepsilon_{1})-
	\mathscr{B}(\bx,t,\bvarepsilon_{2})\|\leq L_{\mathscr{B}}
	\|\bvarepsilon_{1}-\bvarepsilon_{2}\|$
	for all \\ $\bvarepsilon_{1}$, $\bvarepsilon_{2}\in \mathbb{S}^{d}$, a.e. $(\bx, t) \in Q$.
	\smallskip
	\item[(3)]
	${\mathscr B}(\cdot, \cdot, {\bf{0}}) \in L^2(Q; {\mathbb{S}}^d)$.
}

\smallskip

\noindent
$\underline{H({\mathscr G})}:$ \quad
$\mathscr {G}\colon Q \times \mathbb{S}^{d} 
\times \mathbb{S}^d \rightarrow \mathbb{S}^{d}$
is such that

\smallskip

\lista{
	\item[(1)]
	$\mathscr{G}(\cdot,\cdot,\bsigma,\bvarepsilon)$
	is measurable on $Q$ for all
	$\bsigma$, $\bvarepsilon \in \mathbb{S}^{d}$.
	\smallskip
	\item[(2)]
	there exists $L_{\mathscr{G}}>0$ such that 
	for all $\bsigma_1$, $\bsigma_2$, $\bvarepsilon_{1}$, $\bvarepsilon_{2}\in \mathbb{S}^{d}$, 
	a.e. $(\bx, t) \in Q$, it holds \\ 
	$\|\mathscr{G}(\bx,t, \bsigma_1,\bvarepsilon_{1})-
	\mathscr{G}(\bx,t,\bsigma_2, \bvarepsilon_{2})\|
	\le L_{\mathscr{G}}\,
	( \| \bsigma_1 - \bsigma_2\| +
	\|\bvarepsilon_{1}-\bvarepsilon_{2}\|)$.
	\smallskip
	\item[(3)]
	${\mathscr G}(\cdot, \cdot, {\bf{0}}, \bzero) \in L^2(Q; {\mathbb{S}}^d)$.
}

\smallskip

\noindent
$\underline{H({k})}:$ \quad
$k \colon \Gamma_C \times \real \to \real$
is such that

\smallskip

\lista{
	\item[(1)]
	$k(\cdot, r)$ is measurable on
	$\Gamma_{C}$ for all $r \in \real$.
	\smallskip
	\item[(2)]
	there exist $k_1$, $k_2$ such that
	$0 < k_1 \le k(\bx,t)\le k_2$ for all
	$r \in \mathbb{R}$, a.e. $\bx \in \Gamma_{C}$.
	\smallskip
	\item[(3)]
	there is $L_{k}>0$ such that  $|k(\bx,r_{1})-k(\bx,r_{2})|\le L_{k}|r_{1}-r_{2}|$
	for all $r_{1}$, $r_{2}\in \mathbb{R}$, a.e.
	$\bx \in \Gamma_{C}$.
}

\smallskip

\noindent
$\underline{H({j_\nu})}:$ \quad
$j_\nu \colon \Gamma_C  \times  \real \to \real$
is such that

\smallskip

\lista{
	\item[(1)]
	$j_{\nu}(\cdot,r)$ is measurable on
	$\Gamma_{C}$ for all $r \in \mathbb{R}$,
	and $j_{\nu}(\cdot,{e}(\cdot)) \in L^{1}(\Gamma_{C})$
	for some ${e}\in L^{2}(\Gamma_{C})$.
	\smallskip
	\item[(2)]
	$j_{\nu}(\bx,\cdot)$ is locally Lipschitz on
	$\mathbb{R}$ for a.e. $\bx \in \Gamma_{C}$.
	\smallskip
	\item[(3)]
	there is $c_{0} \ge 0$
	such that $|\partial j_{\nu}(\bx,r)|\le
	c_{0}$
	for all $r \in \mathbb{R}$, a.e. $\bx \in \Gamma_{C}$.
	\smallskip
	\item[(4)]
	there exists $\alpha_{j_{\nu}} \ge 0$ such that
	$j_{\nu}^0(\bx,r_{1}; r_{2}-r_{1}) + j_{\nu}^0(\bx,r_{2};r_{1}-r_{2})
	\le \alpha_{j_{\nu}}|r_{1}-r_{2}|^{2}$
	for all $r_{1}$, $r_{2}\in \mathbb{R}$,
	a.e. $\bx \in \Gamma_{C}$.
}

\smallskip

\noindent
$\underline{(H_4)}:$ \quad
$\fb_0 \in L^2(0,T;L^2(\Omega;\real^d))$,
$\fb_N \in L^2(0,T;L^2(\Gamma_N;\real^d))$,  $\bu_0$, $\bw_0\in V$.

\medskip

We introduce the set $K$ of admissible velocity fields 
defined by
\begin{equation}\label{SETK}
K=\{\, \bv \in V \mid
v_{\nu} \leq g \ \ \mbox{\rm on} \ \ \Gamma_{C}\, \},
\end{equation}
and an element $\fb  \in L^2(0,T;V^*)$ given by
\begin{equation}\label{fXXX}
\langle \fb (t),\bv\rangle =
\langle \fb_{0} (t), \bv\rangle_{L^{2}(\Omega;\mathbb{R}^{d})}
+\langle \fb_{N}(t),\bv \rangle_{L^{2}(\Gamma_{N};\mathbb{R}^{d})}
\end{equation}
for all $\bv \in V$, a.e. $t \in (0, T)$.

We now turn to the weak formulation of
Problem~\ref{Plastic}. 
Let $(\bu, \bsigma)$ be a smooth solution to this problem which means that the data are smooth functions such that all the derivatives and all the conditions are satisfied in the usual sense at each point.
Let $\bv\in K$ and $t\in (0,T)$. 
We multiply (\ref{equation2}) by $\bv-\bu'(t)$,
use Green's formula, see~\cite[Theorem~2.25]{MOSBOOK}, 
and apply the boundary conditions (\ref{equation3})
and (\ref{equation4})
to obtain
\begin{eqnarray*}
	&&\hspace{-0.5cm}
	\int_\Omega\bu''(t)\cdot(\bv-\bu'(t)) \,dx
	+\int_\Omega \bsigma(t)\cdot
	\big(
	\bvarepsilon(\bv) -\bvarepsilon(\bu'(t))
	\big)
	\, dx \\ [1mm]
	&&\hspace{-0.5cm} \ \
	=\int_\Omega \fb_0(t)\cdot(\bv-\bu'(t))dx
	+\int_{\Gamma_N} \fb_N(t)\cdot(\bv-\bu'(t))\, d\Gamma
	+\int_{\Gamma_C} \bsigma(t)\bnu \cdot(\bv-\bu'(t))\, d\Gamma.
\end{eqnarray*}
By (\ref{equation5}) and the definition
of the Clarke subgradient, we have
\begin{eqnarray}
&&
\sigma_{\nu}(t) (v_\nu - u_\nu'(t))
= (\sigma_{\nu}(t) + \xi(t))(v_\nu - g)
- (\sigma_{\nu}(t) + \xi(t))(u_\nu'(t) - g)
\nonumber \\[2mm]
&&\qquad
- \xi(t) (v_\nu - u_\nu'(t))
\ge
- k(u_\nu(t)) \, j_\nu^0(u_\nu'(t); v_\nu - u_\nu'(t))
\ \ \mbox{on} \ \ \Gamma_C.
\label{normal}
\end{eqnarray}
On the other hand, (\ref{DEC}) and the frictionless condition (\ref{equation6}) imply
\begin{equation}\label{tangent}
\bsigma (t)\bnu \cdot (\bv - \bu'(t))
=
\sigma_{\nu}(t) (v_\nu - u_\nu'(t))
\ \ \mbox{on} \ \ \Gamma_C.
\end{equation}
Combining (\ref{fXXX}), 
(\ref{normal}) and (\ref{tangent}),
we have
\begin{eqnarray}
&&
\int_\Omega\bu''(t)\cdot(\bv-\bu'(t))\, dx
+\langle \bsigma(t), \bvarepsilon(\bv)-\bvarepsilon(\bu'(t))
\rangle_{\mathcal{H}}
\nonumber \\ [1mm]
&& \qquad\qquad +\int_{\Gamma_{C}}k(u_{\nu}(t))\, j_{\nu}^0(u'_{\nu}(t);
v_{\nu}-u'_{\nu}(t))\, d\Gamma \ge
\langle \fb (t), \bv - \bu'(t) \rangle.
\label{NEW4}
\end{eqnarray}
From the constitutive law (\ref{equation1}), we get 
\begin{equation}\label{ETA1}
\bsigma(t) = {\mathscr A}\bvarepsilon ({\bu}'(t))
+{\mathscr B}\bvarepsilon({\bu}(t)) + \beeta(t)
\end{equation} 
with
\begin{equation}\label{ETA2}
\beeta (t) = \int_0^t{\mathscr G}\big(
\bsigma(s) - {\mathscr A}\bvarepsilon ({\bu}'(s)), 
\bvarepsilon({\bu}(s) \big)\,ds . 
\end{equation}
Inserting (\ref{ETA1}) into (\ref{NEW4}) and using (\ref{ETA2}), we obtain the following variational formulation of Problem~\ref{Plastic}.
\begin{Problem}\label{CONTACTX1}
Find 
$\beeta \colon (0, T) \to {\mathcal{H}}$ and
$\bu \colon (0,T)\to V$ such that 
$\bu'(t) \in K$ a.e. $t\in (0, T)$ 
and
\begin{eqnarray*}
&&\hspace{-0.4cm}
\beeta'(t) = {\mathscr{G}} 
(\beeta(t) + \mathscr{B}(\bvarepsilon(\bu(t))),
\bvarepsilon(\bu(t))) 
\ \ \mbox{\rm for all} \ \ t \in (0, T), 
\\[2mm]
&&\hspace{-0.4cm}
\int_\Omega\bu''(t)\cdot(\bv-\bu'(t)) \, dx
+\langle \mathscr{A}(\bvarepsilon(\bu'(t)))
+ \mathscr{B}(\bvarepsilon(\bu(t))) + \beeta (t),
\bvarepsilon(\bv)-\bvarepsilon(\bu'(t))
\rangle_{\mathcal{H}} \\
&&\hspace{-0.4cm} \quad +\int_{\Gamma_{C}}
k(u_{\nu}(t)) \,
j_{\nu}^0(u'_{\nu}(t); v_{\nu}-u'_{\nu}(t))\, d\Gamma
\ge
\langle \fb(t), \bv - \bu'(t) \rangle
\ \ \mbox{\rm for all} \ \ \bv \in K, 
\ \mbox{\rm a.e.} \ t \in (0,T), \\[1mm]
&&\hspace{-0.4cm}
\beeta (0) = \bzero, \ \bu(0) = \bu_0, \ \bu'(0) = \bw_0.
\end{eqnarray*}
\end{Problem}

Problem~\ref{CONTACTX1} couples the ordinary differential equation with a hemivariational inequality with constraints. 
The result below concerns the unique solvability, 
continuous dependence, and regularity of solution to Problem~\ref{CONTACTX1}.
\begin{Theorem}\label{existence3}
Assume hypotheses
$H({\mathscr A})$, $H({\mathscr B})$,
$H({\mathscr G})$, $H(k)$, $H(j_\nu)$, 
$(H_4)$ and the smallness condition
\begin{equation}\label{CONTACT_smallness}
\alpha_{j_{\nu}}k_2\|\gamma\|^{2} < m_{\mathscr{A}}.
\end{equation}
Then Problem~$\ref{CONTACTX1}$ has a solution 
$\beeta \in H^1(0,T; {\mathcal{H}})$ 
and 
$\bu \in C([0, T]; V)$ such that 
$\bu' \in {\mathbb{W}}$ 
with
$\bu'(t) \in K$ for a.e. $t \in (0, T)$.
If, in addition,
\begin{equation}\label{EXTRA} 
\mbox{\rm either} \ \ j_\nu (\bx, \cdot) 
\ \ \mbox{\rm or} \ \ 
-j_\nu (\bx, \cdot) \ \ \mbox{\rm is regular for a.e.} 
\ \bx \in \Gamma_C,
\end{equation}
then the solution to Problem~$\ref{CONTACTX1}$ 
is unique, and
there exists a constant $c > 0$ such that 
\begin{eqnarray}\label{EST555}
&&\hspace{-0.2cm}
\| \beeta_1 - \beeta_2 \|_{H^1(0, T; E)}
+ \| \bu_1 - \bu_2 \|_{C(0, T; V)}
+ \| \bu'_1 - \bu'_2 \|_{L^2(0, T; V)} \\[1mm]
&&
\hspace{-0.2cm} \quad \le 
c \, (
\| \bu_0 - \widetilde {\bu}_0 \|
+ \| \bw_0 - \widetilde {\bw}_0\|
+ 
\| \fb_0 - {\widetilde{\fb}_0} 
\|_{L^2(0, T; L^2(\Omega;\real^d))}
+ \| \fb_N - {\widetilde{\fb}_N} 
\|_{L^2(0, T; L^2(\Gamma_N; \real^d))}),
\nonumber 
\end{eqnarray}
where $(\beeta_1, \bu_1)$ and $(\beeta_2, \bu_2)$ 
are the unique solutions to Problem~$\ref{CONTACTX1}$ 
corresponding to
\begin{equation*} 
(\bu_0, \bw_0, \fb_0, \fb_N), 
(\widetilde{\bu}_0, {\widetilde{\bw}}_0, 
\widetilde{\fb}, \widetilde{\fb}_N) \in V\times V
\times L^2(0, T; L^2(\Omega;\real^d) \times 
L^2(\Gamma_N; \real^d)),
\end{equation*} 
respectively.
Moreover, we have the regularity $(\beeta, \bu') \in 
C([0,T]; {\mathcal{H}} \times L^2(\Omega;\real^d))$.
%
\end{Theorem}

\noindent
{\bf Proof}. \
Let $\bw(t)=\bu'(t)$ for all $t\in (0,T)$. 
Then $\bu(t) = \bu_0 + \int_0^t \bw(s) \, ds$ for all $t\in (0,T)$ and Problem~\ref{CONTACTX1} 
can be reformulated as follows.
\begin{Problem}\label{CONTACTX33}
Find 
$\beeta \colon (0, T) \to {\mathcal{H}}$ and
$\bw \colon (0,T)\to V$ such that 
$\bw(t) \in K$ a.e. $t\in (0, T)$, 
$\bw(0)=\bw_0$, $\beeta (0) = \bzero$ 
and
\begin{eqnarray*}
&&\hspace{-0.3cm}
\beeta'(t) = {\mathscr{G}} 
\left(
\beeta(t) + \mathscr{B}
\left(
\bvarepsilon(\bu_0 + \int_0^t \bw(s) \, ds)
\right),
\bvarepsilon(\bu_0 + \int_0^t \bw(s) \, ds ) \right) 
\ \ \mbox{\rm for all} \ \ t \in (0, T), 
\\[2mm]
&&\hspace{-0.3cm}
\int_\Omega\bw'(t)\cdot(\bv-\bw(t)) \, dx
+ \langle \mathscr{A}(\bvarepsilon(\bw(t)))
+ \mathscr{B}(\bvarepsilon(\bu_0 + \int_0^t \bw(s) \, ds)) + \beeta (t),
\bvarepsilon(\bv)-\bvarepsilon(\bw(t))
\rangle_{\mathcal{H}} \\
&&\hspace{-0.3cm} \quad +\int_{\Gamma_{C}}
k\left(
u_{0\nu} + \int_0^t w_\nu (s) \, ds \right) \,
j_{\nu}^0(w_{\nu}(t); v_{\nu}-w_{\nu}(t))\, d\Gamma
\ge
\langle \fb (t), \bv - \bw(t) \rangle 
\\[2mm]
&&
\ \ \ \ \mbox{\rm for all} \ \ \bv \in K, 
\ \mbox{\rm a.e.} \ t \in (0,T). 
\end{eqnarray*}
\end{Problem}
Let 
$E := {\mathcal{H}}$, 
$U := V$,
$Z = X := L^2(\Gamma_C)$, 
and  
$A \colon (0,T) \times E \times V \to V^*$,
$S \colon L^2(0,T;V) \to L^2(0,T; U)$,
$R_1 \colon L^2(0,T;V) \to L^2(0,T; V^*)$,
$R_2 \colon L^2(0,T;V) \to L^2(0,T; Z)$, 
$j \colon Z \times X \to \real$, 
$F \colon (0,T) \times E \times U \to E$, 
$M \colon V \to X$ 
be defined by
\begin{eqnarray*}
&&
\langle A (t, \beeta, \bv), \bz \rangle
:= 
\langle\mathscr{A}(t,\bvarepsilon(\bv)) + \beeta,
\bvarepsilon(\bz) \rangle_{\mathcal{H}}
\ \ \mbox{for} \ \ \bv, \bz \in V, \ \beeta \in {\mathcal{H}}, \ \mbox{a.e.} \ t \in (0,T),
\\
&&
(S\bw)(t) := \bu_0 + \int_0^t \bw(s) \, ds 
\ \ \mbox{for} \  \ \bw \in L^2(0,T; V),
\ \mbox{a.e.} \ t \in (0,T),
\\
&&
\langle (R_1\bw)(t), \bv \rangle 
:= \langle {\mathscr{B}} \bvarepsilon 
((S\bw)(t)),
\bvarepsilon (\bv )\rangle_{\mathcal{H}}
\ \ \mbox{for} \  \ \bw \in L^2(0,T; V), \ \bv \in V,
\ \mbox{a.e.} \ t \in (0,T),
\\
&&
(R_2\bw)(t) := ((S\bw)(t))_\nu
= u_{0\nu} + \int_0^t w_\nu(s) \, ds 
\ \ \mbox{for} \  \ \bw \in L^2(0,T; V), 
\ \mbox{a.e.} \ t \in (0,T), 
\\
&&
j(z, v) = \int_{\Gamma_C} k(\bx, z) \,
j_\nu (\bx, v) \, d\Gamma 
\ \ \mbox{for} \ \ z \in Z, \ v \in X, 
\ \mbox{a.e.} \ t \in (0,T),
\\
&&
F(t, \beeta, \bw)(\bx) :=
{\mathscr{G}} (\bx, t, \beeta (\bx) 
+ {\mathscr{B}}(\bvarepsilon(\bw(\bx))), 
\bvarepsilon(\bw(\bx))) 
\ \mbox{for} \ \beeta \in E, \ \bw \in U,
\ \mbox{a.e.} \ t \in (0,T), \\[1mm]
&&
M \bv := v_\nu 
\ \ \mbox{for} \ \ \bv \in V, 
\end{eqnarray*}
and $\varphi \equiv 0$. 
With the above notation,
we consider the following auxiliary system 
associated with Problem~\ref{CONTACTX33}.
\begin{Problem}\label{CONTACTX2}
Find $\beeta \in H^1(0,T; E)$ and 
$\bw\in \mathbb{W}$ such that
$\bw(t) \in K$ for a.e. $t \in (0, T)$, 
$\bw(0)=\bw_0$, $\beeta(0) = \bzero$ 
and
\begin{equation*}
\begin{cases}
\displaystyle
\beeta'(t) = F(t, \beeta(t), (S\bw)(t))
\ \ \mbox{\rm for all} \ \ t \in (0,T),
\\[1mm]
\langle \bw'(t) + A(t, \beeta(t), \bw(t))
+ (R_1 \bw)(t) - \fb(t),
\bv - \bw(t) \rangle
\\[1mm]
\qquad \ \ \,
+ \, j^0 ((R_2\bw)(t), M \bw(t); M\bv
- M\bw(t)) \ge 0
\ \ \mbox{\rm for all} \ 
\bv \in K, \ \mbox{\rm a.e.} \ t \in (0,T).
\end{cases}
\end{equation*}
\end{Problem}

First, we apply Corollary~\ref{Corollary2} to deduce that Problem~\ref{CONTACTX2} is solvable.
We will verify hypotheses
$H(F)$, $H(A)$, $H(j)$, $H(\varphi)$, $H(K)$, $H(M)$, 
and $(H_1)$--$(H_3)$. 
We use assumptions $H(\mathscr{G})$ and 
$H(\mathscr{B})$ to obtain the estimate
\begin{eqnarray*}
&&
\| F(t, \beeta, \bw) \|_{\mathcal{H}}
\le 
\| F(t, \beeta, \bw) - F(t, \bzero, \bzero) \|_{\mathcal{H}}
+ \| F(t, \bzero, \bzero) \|_{\mathcal{H}}
\\[2mm]
&&
\quad \le 
\sqrt{3} \, L_{\mathscr{G}} \,  
(\| \beeta \|_{\mathcal{H}} + L_{\mathscr{B}} \| \bw\| 
+ \| \bw \| ) 
+ \sqrt{2} \, (L_{\mathscr{G}} \, 
\| {\mathscr{B}}(\cdot, t,\bzero) \|_{\mathcal{H}}
+ \| {\mathscr{G}}(\cdot, t,\bzero,\bzero) \|_{\mathcal{H}})
\end{eqnarray*}
for all $\beeta \in {\mathcal{H}}$, 
$\bw \in V$, a.e. $t \in (0, T)$. 
Hence, the function 
$t \mapsto F(t, \beeta, \bw)$ belongs to 
$L^2(0,T;{\mathcal{H}})$ 
for all $\beeta \in {\mathcal{H}}$, $\bw \in V$.
Similarly, we find that 
\begin{equation*}
\| F(t, \beeta_1, \bw_1)-F(t,\beeta_2, \bw_2) \|_{\mathcal{H}}
\le 
\sqrt{2} \, L_{\mathscr{G}} \,  
(\sqrt{2}\, \| \beeta_1-\beeta_2 \|_{\mathcal{H}} + \sqrt{2} \, L_{\mathscr{B}} \| \bw_1-\bw_2\| 
+ \| \bw_1-\bw_2 \| ) 
\end{equation*}
for all $\beeta_1$, $\beeta_2 \in {\mathcal{H}}$, 
$\bw_1$, $\bw_2 \in V$, a.e. $t \in (0, T)$.
We deduce that the map $F$ satisfies $H(F)$.

Next, we show that operator ${\mathscr{A}}$ satisfies $H(A)$. It is clear from $H({\mathscr{A}})$ that 
$A(\cdot, \beeta, \bv)$ is measurable for all 
$\beeta\in E$, $\bv \in V$, 
and 
$A(t,\cdot, \bv)$ is continuous for all $\bv \in V$, 
a.e. $t \in (0, T)$. 
From $H({\mathscr{A}})$(3), we have
$
\| {\mathscr{A}}(t,\bvarepsilon(\bv)) \|_{\mathcal H} 
\le \sqrt{2} \, 
(\| {\widetilde{a}}_0(t)\|_{L^2(\Omega)} + 
{\widetilde{a}}_2 \, 
\| \bvarepsilon (\bv) \|_{\mathcal H})
$
for all $\bv \in V$, a.e. $t \in (0, T)$, 
and consequently
\begin{eqnarray*}
&&
\| A(t, \beeta, \bv) \|_{V^*} 
= \sup_{\| z\|\le 1} \langle
{\mathscr{A}}(t, \bvarepsilon (\bv))+\beeta, \bvarepsilon(\bz) \rangle_{\mathcal H} 
\le 
\sup_{\| z \|\le 1}
\left(
\| {\mathscr{A}}(t, \bvarepsilon (\bv))\|_{\mathcal H} 
+ \| \beeta \|_{\mathcal H} 
\right) \, \| \bz \| \\
&& 
\qquad \le
\sqrt{2} \, 
\| {\widetilde{a}}_0(t)\|_{L^2(\Omega)} + 
\| \beeta \|_{\mathcal H} +
\sqrt{2} \, {\widetilde{a}}_2 \, \| \bv \|, 
\end{eqnarray*}
which implies condition $H(A)$(c)  
with 
$a_0(t) = \sqrt{2} \, 
\| {\widetilde{a}}_0(t)\|_{L^2(\Omega)}$, 
$a_1 = 1$ and 
$a_2 = \sqrt{2} \, {\widetilde{a}}_2$. 
%
To demonstrate $H(A)$(d), we observe that 
\begin{equation*}
\| A(t, \beeta_1, \bv) - A(t, \beeta_2, \bv) \|_{V^*} 
= \sup_{\| z\|\le 1} \langle
\beeta_1 - \beeta_2, \bvarepsilon(\bz) 
\rangle_{\mathcal H} 
\le 
\sup_{\| z \|\le 1}
\| \beeta_1 - \beeta_2\|_{\mathcal H} 
\| \bvarepsilon (\bz) \|_{\mathcal H} 
\le \| \beeta_1 - \beeta_2\|_{\mathcal H}
\end{equation*}
for all $\beeta_1$, $\beeta_2 \in {\mathcal H}$, 
$\bv \in V$, a.e. $t \in (0, T)$, 
and from $H({\mathscr{A}})$(4), we deduce
\begin{eqnarray*}
&&
\langle A(t, \beeta, \bv_1) - A(t, \beeta, \bv_2)
\rangle = 
\langle
{\mathscr{A}}(t, \bvarepsilon (\bv_1))-
{\mathscr{A}}(t, \bvarepsilon (\bv_2)),
\bvarepsilon (\bv_1)-\bvarepsilon (\bv_2)
\rangle_{\mathcal H} \\
&&
\qquad 
\ge m_{\mathscr{A}} \int_{\Omega}
\| \bvarepsilon (\bv_1)-\bvarepsilon (\bv_2)\|^2_{\mathbb{S}^d} \, dx = 
m_{\mathscr{A}} \, \| \bv_1-\bv_2\|^2
\end{eqnarray*}
for all $\beeta \in {\mathcal H}$, $\bv_1$,  
$\bv_2 \in V$, a.e. $t \in (0, T)$.
Exploiting the last two inequalities, by Remark~\ref{REM2}, we obtain $H(A)$(d) with $m_A = m_{\mathscr{A}}$ 
and ${\bar{m}}_A = 0$. Hence $H(A)$ holds. 

From hypotheses $H(k)$, $H(j_\nu)$, 
and~\cite[Theorem~3.47]{MOSBOOK},
it follows that function $j$ satisfies $H(j)$(a)--(d), and
\begin{equation}\label{j-inequality}
j^0(z, v; w) \le \int_{\Gamma_{C}} 
k(\bx, z) \, j_\nu^0 (\bx, v; w) \, d\Gamma 
\ \ \mbox{for all} \ \ z, v, w \in X.
\end{equation}
We use (\ref{j-inequality}) with $H(k)$ and 
$H(j_\nu)$ to get
\begin{eqnarray*}
&&\hspace{-0,5cm}
j^0(z_1, v_1; v_2 - v_1) + j^0(z_2, v_2; v_1 - v_2) \\ [1mm]
&&
\le \int_{\Gamma_C} \Big(
k(z_1) j_\nu^0 (v_{1}; v_{2} - v_{1}) + 
k(z_2) j_\nu^0 (v_{2}; v_{1} - v_{2}) \Big) 
\, d\Gamma  \\ [1mm]
&&
\le 
\int_{\Gamma_C} 
\big( k(z_1) - k(z_2) \big) 
j_\nu^0 (v_{1}; v_{2} - v_{1}) + 
\, k(z_2) \Big( 
j_\nu^0 (v_{1}; v_{2} - v_{1}) +
j_\nu^0 (v_{2}; v_{1} - v_{2}) \Big) \, d\Gamma \\ [1mm]
&&
\le
c_0 \, L_k \int_{\Gamma_C} 
| z_1(\bx) - z_2(\bx)| \, \| v_1(\bx) - v_2(\bx)\| \, d\Gamma + 
\alpha_{j_{\nu}} \, k_2 
\int_{\Gamma_C} \| v_1(\bx) - v_2(\bx) \|^2 \, d\Gamma 
\\ [1mm]
&&
\le 
{\bar{m}}_{j} \| z_1 - z_2 \|_Z \| v_1 - v_2 \|_X 
+ m_{j} \| v_1 - v_2 \|^2_X,
\end{eqnarray*}
where 
$m_{j} = \alpha_{j_{\nu}} k_2$ and
${\bar{m}}_{j} = c_0 L_k$. 
Hence $H(j)$(e) follows, and therefore,  
$j$ satisfies $H(j)$.

Condition $H(\varphi)$ holds trivially. 
It is clear that the set $K$ is a closed and convex subset of $V$ with $\bzero \in K$, i.e., $H(K)$ holds.
By the linearity and continuity of the normal trace operator, it is obvious that $M_1 = M$ and $H(M)$(a) is satisfied. 
Condition $H(M)$(b) holds as a consequence of the proof 
of~\cite[Theorem~2.18, p.59]{AMMA2}. 
Using the regularity of $\fb_0$ and $\fb_N$ in $(H_4)$, 
we readily obtain that the functional defined by (\ref{fXXX}) satisfies $\fb \in L^2(0, T; V^*)$.

Condition $(H_1)$ is a consequence of the smallness 
hypothesis (\ref{CONTACT_smallness}) while 
condition $(H_2)$ follows from $(H_4)$.
Finally, 
we argue as in~\cite[Theorem~14.2, p.367]{AMMA14}
to demonstrate that operators $S$, $R_1$ and $R_2$ satisfy the inequalities $(H_3)$(a)--(c), respectively.

Having verified all hypotheses 
of Corollary~\ref{Corollary2}, we deduce from it that Problem~\ref{CONTACTX2} has a unique solution
$\beeta \in H^1(0,T; E)$, $\bw\in \mathbb{W}$ 
such that $\bw(t) \in K$ for a.e. $t \in (0, T)$.
By the inequality (\ref{j-inequality}) it follows that any solution to Problem~\ref{CONTACTX2} is a solution to 
Problem~\ref{CONTACTX33}. 
Moreover, from the relation
$\bu(t) = \int_0^t \bw(s) \, ds + \bu_0$ 
for all $t \in (0, T)$, we conclude that
$\bu \in C([0, T]; V)$, $\bu' \in \mathbb{W}$ with
$\bu'(t) \in K$ for a.e. $t \in (0, T)$.
The proof of existence of solution to Problem~\ref{CONTACTX1} is completed.
 
Next, we suppose the regularity hypothesis 
(\ref{EXTRA}). 
Invoking~\cite[Theorem~3.47(vii)]{MOSBOOK}, 
we deduce that (\ref{j-inequality}) holds 
with equality. This implies that
Problems~\ref{CONTACTX33} and~\ref{CONTACTX2} 
are equivalent. Therefore, in this case,  Problem~\ref{CONTACTX1} is uniquely solvable.

Let $(\beeta_1, \bu_1)$ and $(\beeta_2, \bu_2)$ 
be the unique solutions to Problem~$\ref{CONTACTX1}$ 
corresponding to the data
\begin{equation*} 
(\bu_0, \bw_0, \fb_0, \fb_N), 
(\widetilde{\bu}_0, {\widetilde{\bw}}_0, 
\widetilde{\fb}, \widetilde{\fb}_N) \in V\times V
\times L^2(0, T; L^2(\Omega;\real^d) \times 
L^2(\Gamma_N; \real^d)),
\end{equation*} 
respectively. 
From (\ref{fXXX}) and $(H_4)$, we readily have 
\begin{equation}\label{*3}
\| \fb - \widetilde{\fb} \|_{L^2(0, T; V^*)} 
\le c \, 
\left(
\| \fb_0 - {\widetilde{\fb}_0} 
\|_{L^2(0, T; L^2(\Omega;\real^d))}
+ \| \fb_N - {\widetilde{\fb}_N} 
\|_{L^2(0, T; L^2(\Gamma_N; \real^d))}
\right)
\end{equation}
with a constant $c > 0$, which by (\ref{EST444}) of Corollary~\ref{Corollary2} implies
\begin{eqnarray}\label{*4}
&&\hspace{-0.2cm}
\| \beeta_1 - \beeta_2 \|_{H^1(0, T; E)}
+ \| \bw_1 - \bw_2 \|_{L^2(0, T; V)}
\\[1mm]
&&
\hspace{-0.2cm} \quad \le 
c \, (
\| \bw_0 - \widetilde {\bw}_0\|
+ 
\| \fb_0 - {\widetilde{\fb}_0} 
\|_{L^2(0, T; L^2(\Omega;\real^d))}
+ \| \fb_N - {\widetilde{\fb}_N} 
\|_{L^2(0, T; L^2(\Gamma_N; \real^d))}).
\nonumber 
\end{eqnarray}
Using again the equation 
$\bu(t) = \bu_0 + \int_0^t \bw(s) \, ds$ 
for all $t\in (0,T)$, we obtain 
\begin{equation}\label{*55}
\| \bu_1 - \bu_2 \|_{C(0, T; V)} 
\le \| \bu_0 - \widetilde{\bu}_0 \| 
+ \sqrt{T} \, \| \bw_1 - \bw_2 \|_{L^2(0, T; V)}.
\end{equation}
Adding (\ref{*4}) and (\ref{*55}), and combining with 
$\bw_i = \bu'_i$ for $i=1$, $2$, we deduce 
the inequality (\ref{EST555}).
Finally, we use the continuity of the embeddings 
$$
H^1(0, T; {\mathcal{H}}) \subset C([0, T]; {\mathcal{H}}) 
\ \ \mbox{and} \ \ {\mathbb{W}} \subset 
C([0,T]; L^2(\Omega;\real^d))
$$
to conclude the desired regularity of the solution 
$(\beeta, \bu') \in 
C([0,T]; {\mathcal{H}} \times L^2(\Omega;\real^d))$. 
This completes the proof.
\hfill$\Box$

%


\section{A dynamic frictional viscoelastic contact problem}\label{Application2}

In this section we consider the dynamic contact problem
for viscoelastic materials with friction and adhesion. 
The physical framework is similar to the one considered in the previous section. 
The difference arises in the fact that
now the body $\Omega \subset \real^d$ 
is assumed to be viscoelastic with long memory, 
the process is frictional and it is governed by a convex potential, and the adhesion process is modeled by 
an ordinary differential equation for the bonding variable on a contact surface $\Gamma_C$.
The classical formulation of the contact problem we study in this section is the following.
\begin{Problem}\label{VISCO2}
Find a displacement field
$\bu \colon Q \to\mathbb{R}^d$ and a stress field
$\bsigma \colon Q \rightarrow \mathbb{S}^d$ and 
a bonding field $\beta \colon \Gamma_C \times (0, T) 
\to [0,1]$
such that for all $t\in (0,T)$,
\begin{align}
\label{equation1x}
\bsigma(t) &={\mathscr A}\bvarepsilon ({\bu}'(t))
+{\mathscr B}\bvarepsilon({\bu}(t))
+ \int_0^t {\mathscr{C}}(t-s) \bvarepsilon (\bu'(s))\, ds 
\quad
&{\rm in}\
&\Omega,\\[1mm]
\label{equation2x}
{\bu}''(t)&={\rm Div}\,\bsigma(t)+\fb_0(t)\quad&{\rm in}\ &\Omega,\\[1mm]
\label{equation3x} \bu(t)&=\bzero &{\rm on}\ &\Gamma_D,\\[1mm]
\label{equation4x} \bsigma(t)\bnu&=\fb_N(t)\quad&{\rm on}\ &\Gamma_N,\\[1mm]
\label{equation5x}
u'_{\nu}(t)&\le g, \
\sigma_{\nu}(t) \le 0, \ (u'_{\nu}(t)-g) \,  \sigma_{\nu}(t)=0 
\quad&{\rm	on}\
&\Gamma_C,\\[1mm]
\label{equation6x} 
-\bsigma_\tau(t) &\in \mu(u_\nu(t)) \, h_1(\beta(t)) \,
\partial h_2 (\bu'_\tau(t))
\quad&{\rm on}\ &\Gamma_C,\\[1mm]
\label{equation7x}
\beta'(t)&= G(t, \beta(t), \bu(t))
\quad&{\rm	on}\
&\Gamma_C,\\[1mm]
\label{equation8x}
\beta(0)&= \beta_0 \quad&{\rm	on}\
&\Gamma_C,\\[1mm]
\label{equation9x}
\bu(0)&= \bu_0, \ \ \bu'(t) = \bw_0 
\quad&{\rm	in}\
&\Omega.
%
\end{align} 
\end{Problem}


\noindent
In this problem, the viscoelastic constitutive law 
is of the form (\ref{equation1x}), where 
$\mathscr{A}$, $\mathscr{B}$ and $\mathscr C$ 
denote the viscosity operator, the elasticity operator,
and the relaxation tensor, respectively.

In the study of Problem~\ref{VISCO2} we suppose that the viscosity operator ${\mathscr{A}}$, the elasticity operator ${\mathscr{B}}$ satisfy condition 
$H({\mathscr{A}})$ and $H({\mathscr{B}})$, respectively.
We also assume that the densities of the body 
forces ${\fb_0}$ and tractions ${\fb_N}$, and the initial data $\bu_0$ and $\bw_0$ satisfy hypothesis $(H_4)$. 
The condition (\ref{equation6x}) represents the friction condition which is a law between the tangential velocity and the tangential stress corresponding to the 
friction force. 
This law involves the product of a function $h_1$ 
which depends on the adhesion field $\beta$ and the subdifferential $\partial h_2$ of a convex potential 
$h_2 (\bx, \cdot)$. 
The coefficient $\mu$ is a given function which depends on the normal displacement.
Since condition (\ref{equation5x}) is the Signorini  unilateral contact boundary condition for the normal velocity with a constant $g > 0$, 
analogously as in the previous section, 
we need the set $K$ of admissible unilateral velocity constraints defined by (\ref{SETK}).

The unknown surface variable $\beta$ is called the adhesion field, see~\cite{SHS}. It describes intensity of
adhesion on the contact surface and is governed 
by an ordinary differential equation (\ref{equation7x}) 
with an initial condition (\ref{equation8x}). The value 
$\beta = 1$ means that the adhesion is complete, 
$\beta = 0$ means that there is no adhesion,
and $0<\beta<1$ means that the adhesion is partial.
Further, we need the following hypotheses on the relaxation tensor ${\mathscr{C}}$, 
the potential $h$,
the coefficient $\mu$,
and the adhesive evolution rate function $G$.

\smallskip

\noindent
$\underline{H({\mathscr C})}:$ \quad
$\mathscr {C}\colon Q \times \mathbb{S}^{d}
\rightarrow \mathbb{S}^{d}$
is such that

\smallskip

\lista{
\item[(1)]
$\mathscr{C}(\bx, t,\bvarepsilon) = (c_{ijkl}(\bx,t)\varepsilon_{kl})$
for all	$\bvarepsilon = (\varepsilon_{ij}) \in \mathbb{S}^{d}$, a.e. $(\bx, t) \in Q$.
\smallskip
\item[(2)]
$c_{ijkl}(\cdot, t) = c_{jikl}(\cdot, t)
= c_{lkij}(\cdot, t) \in L^\infty(\Omega)$ 
for a.e. $t \in (0, T)$.
\smallskip
\item[(3)]
$t \mapsto c_{ijkl}(\cdot, t) \in L^\infty(0,T; L^\infty(\Omega))$.
}

\smallskip

\noindent
$\underline{H(h_1)}:$ \quad
$h_1 \colon \Gamma_C \times \real \rightarrow \real$
is such that

\smallskip

\lista{
\item[(1)]
$h_1(\cdot, r)$ is measurable on $\Gamma_C$, for all 
$r \in \real$.
\smallskip
\item[(2)]
$h_1(\bx, \cdot)$ is Lipschitz continuous on $\real$, 
a.e. $\bx \in \Gamma_C$.
\smallskip
\item[(3)]
$0 \le h_1(\bx, r) \le h_0$ for all $r \in \real$, a.e. $\bx \in \Gamma_C$ with $h_0 > 0$.
}

\smallskip

\noindent
$\underline{H(h_2)}:$ \quad
$h_2 \colon \Gamma_C \times \real^{d}
\rightarrow \real$
is such that

\smallskip

\lista{
	\item[(1)]
	$h_2(\cdot, \bxi)$ is measurable on $\Gamma_C$ for 
	$\bxi \in \real^d$.
	\smallskip
	\item[(2)]
	$h_2(\bx, \cdot)$ is convex, 
	lower semicontinuous, and Lipschitz, 
	for a.e. $\bx \in \Gamma_C$.
	\smallskip
	\item[(3)]
	$\|\partial h_2(\bx, \bxi)\| \le 
	\rho_1(\bx) + \rho_2 \| \bxi \|$ for all 
	$\bxi \in \real^d$, a.e. $\bx \in \Gamma_C$ with 
	$\rho_1 \in L^2(\Gamma_C)$, $\rho_1$, $\rho_2 \ge 0$.
	}

\smallskip

\noindent
$\underline{H(\mu)}:$ \quad
$\mu \colon \Gamma_C \times \real \rightarrow \real$
is such that

\smallskip

\lista{
\item[(1)]
$\mu(\cdot, r)$ is measurable on $\Gamma_C$ for all $r \in \real$.
\smallskip
\item[(2)]
$\mu(\bx, \cdot)$ is continuous on $\real$ 
for a.e. $\bx \in \Gamma_C$.
\smallskip
\item[(3)]
$0 \le \mu(\bx, r) \le \mu_0$ for all $r \in \real$,
a.e. $\bx \in \Gamma_C$ with $\mu_0 > 0$.
}

\smallskip

\noindent
$\underline{H(G)}:$ \quad
$G \colon \Gamma_C \times (0, T) \times \real \times \real^d 
\rightarrow \real$
is such that

\smallskip

\lista{
\item[(1)]
$G(\cdot, \cdot, r, \bxi)$ is measurable on 
$\Gamma_C \times (0, T)$ for all $r \in \real$, $\bxi \in \real^d$.
\smallskip
\item[(2)]
$|G(\bx, t, r_1, \xi_1) - G(\bx, t, r_2, \bxi_2)|
\le L_G \, \left( |r_1 - r_2| + \| \bxi_1 - \bxi_2\|\right)$
for all $(r_i, \bxi_i) \in \real \times \real^d$, 
$i=1$, $2$, a.e. $(\bx, t) \in \Gamma_C \times (0, T)$ 
with $L_G > 0$.
\smallskip
\item[(3)]
$G(\bx, t, 0, \bxi) = 0$, 
$G(\bx, t, r, \bxi) \ge 0$ for $r \le 0$, and
$G(\bx, t, r, \bxi) \le 0$ for $r \ge 1$, 
all $\bxi \in \real^d$, 
a.e. $(\bx, t) \in \Gamma_C \times (0, T)$.
}

\smallskip

Let $V$ and ${\mathcal H}$ be the spaces defined by (\ref{SPACESVH}), and $\fb  \in L^2(0,T;V^*)$ be given 
by (\ref{fXXX}).
We use the standard procedure as in Section~\ref{Application1} to obtain the following weak  formulation of Problem~\ref{VISCO2}.

\begin{Problem}\label{CONTACTX88}
	Find 
	$\beta \colon (0, T) \to [0, 1]$ and
	$\bu \colon (0,T)\to V$ such that 
	$\bu'(t) \in K$ a.e. $t\in (0, T)$ 
	and
\begin{eqnarray*}
&&\hspace{-0.4cm}
\beta'(t) = G(t, \beta(t), \bu(t)) 
\ \ \mbox{\rm on} \ \ \Gamma_C 
\ \mbox{\rm for all} \ \, t \in (0, T),  
\\[1mm]
&&\hspace{-0.4cm} 
\langle \bu''(t), \bv-\bu'(t) \rangle 
+\langle \mathscr{A}(\bvarepsilon(\bu'(t)))
+ \mathscr{B}(\bvarepsilon(\bu(t))) 
+ \int_0^t {\mathscr{C}}(t-s) 
\bvarepsilon (\bu'(s)) \, ds,
\bvarepsilon(\bv)-\bvarepsilon(\bu'(t))
\rangle_{\mathcal{H}}
\\
&&\hspace{-0.4cm} \quad
+\int_{\Gamma_{C}}
\mu(u_{\nu}(t)) \, h_1(\beta(t)) \,
\left(
h_2(\bv_\tau)-h_2(\bu'_\tau(t))
\right)	\, d\Gamma 
\ge
\langle \fb(t), \bv - \bu'(t) \rangle
\\[2mm]
&&\hspace{-0.4cm} \qquad
\ \ \mbox{\rm for all} \ \, \bv \in K, 
\ \mbox{\rm a.e.} \ t \in (0,T), \\[1mm]
&&\hspace{-0.4cm}
\beta (0) = \beta_0, \ \bu(0) = \bu_0, \ \bu'(0) = \bw_0.
	\end{eqnarray*}
\end{Problem}

Problem~\ref{CONTACTX88} consists with the ordinary differential equation on the contact surface and a dynamic variational inequality with constraints. 
The following result deals with the unique solvability, 
continuous dependence, and regularity of solution to Problem~\ref{CONTACTX88}.
\begin{Theorem}\label{existence88}
Assume hypotheses
$H({\mathscr A})$, $H({\mathscr B})$,
$H({\mathscr C})$, $H(h_1)$, $H(h_2)$, $H(\mu)$, $H(G)$,  
$(H_4)$ and $\beta_0 \in L^2(\Gamma_C)$ with
$0 \le \beta_0(\bx) \le 1$ for a.e. $\bx \in \Gamma_C$. 
Then Problem~$\ref{CONTACTX88}$ has a unique 
solution 
	$\beta \in H^1(0,T; L^2(\Gamma_C))$ 
	and 
	$\bu \in C([0, T]; V)$ such that 
	$0 \le \beta(\bx, t) \le 1$ for all $t \in (0, T)$, 
	a.e. $\bx \in \Gamma_C$, 
	$\bu' \in {\mathbb{W}}$ 
	with
	$\bu'(t) \in K$ for a.e. $t \in (0, T)$.
Moreover, let  
$(\beta_1, \bu_1)$ and $(\beta_2, \bu_2)$ be the unique solutions to Problem~$\ref{CONTACTX88}$ 
corresponding to the data
\begin{equation*} 
(\beta_0, \bu_0, \bw_0, \fb_0, \fb_N), 
(\widetilde{\beta}_0, \widetilde{\bu}_0, {\widetilde{\bw}}_0, 
\widetilde{\fb}, \widetilde{\fb}_N) \in L^2(\Gamma_C) 
\times V\times V \times L^2(0, T; L^2(\Omega;\real^d) 
\times L^2(\Gamma_N; \real^d)),
\end{equation*}
respectively. Then, there is a constant $c > 0$ such that
\begin{eqnarray*}\label{EST789}
&&\hspace{-0.6cm}
\| \beta_1 - \beta_2 \|_{H^1(0, T;L^2(\Gamma_C))}
+ \| \bu_1 - \bu_2 \|_{C(0, T; V)}
+ \| \bu'_1 - \bu'_2 \|_{L^2(0, T; V)} \le 
c \, ( \| \beta_0 - \widetilde{\beta}_0 \|_{L^2(\Gamma_C)}
\\[1mm]
&&
\hspace{-0.4cm} \quad  
+ \, \| \bu_0 - \widetilde {\bu}_0 \|
+ \| \bw_0 - \widetilde {\bw}_0\|
+ 
\| \fb_0 - {\widetilde{\fb}_0} 
\|_{L^2(0, T; L^2(\Omega;\real^d))}
+ \| \fb_N - {\widetilde{\fb}_N} 
\|_{L^2(0, T; L^2(\Gamma_N; \real^d))}),
\nonumber 
\end{eqnarray*}
Moreover, $(\beta, \bu') \in 
C([0,T]; L^2(\Gamma_C) \times L^2(\Omega;\real^d))$.
\end{Theorem}
\noindent
{\bf Proof}. \ 
Let $\bw(t)=\bu'(t)$ for all $t\in (0,T)$, hence 
$\bu(t) = \bu_0 + \int_0^t \bw(s) \, ds$ 
for all $t\in (0,T)$.
We reformulate Problem~$\ref{CONTACTX88}$ in terms of the velocity $\bw$ using the notation below.

Let 
$E := Y = L^2(\Gamma_C)$, 
$U := V$, $X = L^2(\Gamma_C;\real^d)$,
$A \colon (0,T) \times V \to V^*$,
$I \colon L^2(0, T; V) \to L^2(0, T; V)$,
$S \colon L^2(0,T;V) \to L^2(0,T; U)$,
$R_1 \colon L^2(0,T;V) \to L^2(0,T; V^*)$,
$R_3 \colon L^2(0,T;V) \to L^2(0,T; Y)$,
$\varphi \colon (0, T) \times E \times Y \times X \to \real$,
$F \colon (0,T) \times E \times U \to E$, 
$M \colon V \to X$ 
be defined by
\begin{eqnarray*}
&&
\langle A (t, \bv), \bz \rangle := 
\langle\mathscr{A}(t,\bvarepsilon(\bv)), 
\bvarepsilon(\bz) \rangle_{\mathcal{H}}
\ \ \mbox{for} \ \ \bv, \bz \in V, 
\ \mbox{a.e.} \ t \in (0,T),
\\
&&
(I\bw)(t) := \bu_0 + \int_0^t \bw(s) \, ds 
\ \ \mbox{for} \  \ \bw \in L^2(0,T; V),
\ \mbox{a.e.} \ t \in (0,T), 
\\
&&
(S\bw)(t) := ((I\bw)(t))_\tau 
= \bu_{0\tau} + \int_0^t \bw_\tau(s) \, ds 
\ \ \mbox{for} \  \ \bw \in L^2(0,T; V),
\ \mbox{a.e.} \ t \in (0,T),
\\
&&
\langle (R_1\bw)(t), \bv \rangle 
:= \langle {\mathscr{B}} \bvarepsilon ((I\bw)(t)) 
+ \int_0^t {\mathscr{C}}(t-s) \bvarepsilon (\bw(s)) 
\, ds, \bvarepsilon (\bv )\rangle_{\mathcal{H}}
\\
&&\qquad\qquad\qquad\qquad 
\ \ \mbox{for} \  \ \bw \in L^2(0,T; V), \ \bv \in V,
\ \mbox{a.e.} \ t \in (0,T),
\\
&&
(R_3\bw)(t) := ((I\bw)(t))_\nu
= u_{0\nu} + \int_0^t w_\nu(s) \, ds 
\ \ \mbox{for} \  \ \bw \in L^2(0,T; V), 
\ \mbox{a.e.} \ t \in (0,T), 
\\
&&
\varphi(t, \beta, y, \bv)
:= \int_{\Gamma_{C}} \mu (y) \, h_1(\beta) \, 
h_2 (\bv_\tau) \, d\Gamma 
\ \ \mbox{for} \ \ \beta \in E, \ y \in Y, \ \bv \in X, 
\ \mbox{a.e.} \ t \in (0,T),
\\
&&
F(t, \beta, \bw)(\bx) :=
G(\bx, t, \beta (\bx, t), \bw_\tau (\bx,t))
\ \mbox{for} \ \beta \in E, \ \bw \in U,
\ \mbox{a.e.} \ t \in (0,T), \\[1mm]
&&
M \bv := \bv_\tau 
\ \ \mbox{for} \ \ \bv \in V, 
\end{eqnarray*}
and $j \equiv 0$. 
Under this notation Problem~$\ref{CONTACTX88}$ 
is equivalently written as
\begin{Problem}\label{CONTACTX88a}
Find $\beta \in H^1(0,T; E)$ and 
	$\bw\in \mathbb{W}$ such that
	$\bw(t) \in K$ for a.e. $t \in (0, T)$, 
	$\bw(0)=\bw_0$, $\beta(0) = \beta_0$ 
	and
	\begin{equation*}
	\begin{cases}
	\displaystyle
	\beta'(t) = F(t, \beta(t), (S\bw)(t))
	\ \ \mbox{\rm for all} \ \ t \in (0,T),
	\\[1mm]
\langle \bw'(t) + A(t, \bw(t))
+ (R_1 \bw)(t) - \fb(t),
\bv - \bw(t) \rangle 
+ \, \varphi(t, \beta(t), (R_3\bw)(t), M \bv)
\\[1mm]
\qquad \ \ \,
- \, \varphi(t, \beta(t), (R_3\bw)(t), M \bw(t)) \ge 0
	\ \ \mbox{\rm for all} \ 
	\bv \in K, \ \mbox{\rm a.e.} \ t \in (0,T).
	\end{cases}
	\end{equation*}
\end{Problem}
Now we apply Corollary~\ref{Corollary2} to deduce the well-posedness of Problem~\ref{CONTACTX88a}.  
To this end, we are able to verify hypotheses
$H(F)$, $H(A)$, $H(j)$, $H(\varphi)$, $H(K)$, $H(M)$, 
$\fb \in L^2(0, T; V^*)$, and $(H_1)$--$(H_3)$. 
Note that hypotheses $H(j)$ and $(H_1)$ hold trivially and no smallness condition is needed. 
Condition $H(F)$ is a consequence of $H(G)$ and can be proved analogously as in the proof of 
Theorem~\ref{existence3}. 
Analogously, we can also show that $H(A)$, $H(K)$, $H(M)$, and $(H_3)$ are satisfied, while $(H_2)$ follows from $(H_4)$.
Additionally, a careful examination 
of~\cite[Lemma~5]{MO2008} implies that 
$0 \le \beta(\bx, t) \le 1$ for all $t \in (0, T)$, 
a.e. $\bx \in \Gamma_C$. 
Next, we use 
$H(h_1)$, $H(h_2)$ and $H(\mu)$ to check that  $H(\varphi)$(a)--(d) holds with 
$c_{0\varphi} = h_0 \mu_0 \| \rho\|_{L^2(\Gamma_C)}$,
$c_{1\varphi} = c_{2\varphi} = 0$ and 
$c_{3\varphi} = h_0 \mu_0 \rho_2$. 
Condition $H(\varphi)$(e) follows from the inequality 
\begin{eqnarray*}
&&\hspace{-0.4cm}
\varphi (t, \beta_1, y_1, v_2) - 
\varphi (t, \beta_1, y_1, v_1) + 
\varphi (t, \beta_2, y_2, v_1) - 
\varphi (t, \beta_2, y_2, v_2) \\
&&
\le 
\mu_0 \int_{\Gamma_C} 
|h_1(\beta_1)-h_1(\beta_2)|
\, | h_2(\bv_{1\tau}) - h_2(\bv_{2\tau}) | \, d\Gamma
\le 
c \, \| \beta_1 - \beta_2\|_{E} \, 
\| \bv_1 - \bv_2\|_X
\end{eqnarray*}
for all $\beta_1$, $\beta_2 \in E$, 
$\bv_1$, $\bv_2 \in X$ with $c > 0$. 
Hence hypothesis $H(\varphi)$ is verified.
By Corollary~\ref{Corollary2} we deduce the unique solvability of Problem~\ref{CONTACTX88a}. 

Moreover, due to (\ref{EST444}) in Corollary~\ref{Corollary2}, we have
\begin{eqnarray}\label{**4}
&&\hspace{-0.2cm}
\| \beta_1 - \beta_2 \|_{H^1(0, T;L^2(\Gamma_C))}
+ \| \bw_1 - \bw_2 \|_{L^2(0, T; V)} \le 
c \, (\| \beta_0 - \widetilde{\beta}_0 \|_{L^2(\Gamma_C)}
\\[1mm]
&&
\hspace{-0.2cm} \quad 
+ \, \| \bw_0 - \widetilde {\bw}_0\|
+ \| \fb_0 - {\widetilde{\fb}_0} 
\|_{L^2(0, T; L^2(\Omega;\real^d))}
+ \| \fb_N - {\widetilde{\fb}_N} 
\|_{L^2(0, T; L^2(\Gamma_N; \real^d))}).
\nonumber 
\end{eqnarray}
We combine (\ref{*55}) with (\ref{**4}) to deduce 
the inequality in the thesis of the theorem. 
As in Theorem~\ref{existence3}, we use the 
continuous embeddings 
$H^1(0, T; L^2(\Gamma_C)) \subset C([0, T]; L^2(\Gamma_C))$  and ${\mathbb{W}} \subset C([0,T]; L^2(\Omega;\real^d))$
to conclude the desired regularity of the solution 
$(\beta, \bu') \in 
C([0,T]; L^2(\Gamma_C) \times L^2(\Omega;\real^d))$. 
This completes the proof.
\hfill$\Box$

\medskip

We refer to~\cite[Chapter~11.4]{SST}, 
\cite[Chapter~5]{SHS}, and 
\cite{Bartosz,HLM2015,MO2008,MShengda2018} 
for examples of the adhesive evolution rate function 
$G$ which depends on both the bonding field $\beta$ and
the displacement and may change sign. This allows for rebonding to take place after debonding, and it allows for
possible cycles of debonding and rebonding.
Note also that the choice 
$h_2(\bxi) = \| \bxi\|$ for $\bxi \in \real^d$ 
leads to a modified version of Coulomb's law which is usually used to model the frictional contact.

\section{Final comments}

We note that similar well-posedness results for more 
complicated contact models can be obtained when various conditions on different parts of the contact boundary 
are considered.

Furthermore, by the application of Theorem~\ref{Theorem1} and Corollary~\ref{Corollary2}, 
one can deal with the well-posedness of dynamic frictional and frictionless contact models involving the internal state variables in viscoplasticity, see~\cite[Chapter~3]{SOFMIG}, 
and the wear phenomena in contact problems, 
see e.g.~\cite[Chapter~11]{SST}, \cite[Chapter~2]{SHS}.

It would be interesting to address the following open issues related to the results of this paper. First it is of interest
to  examine the penalty methods  
for the system (\ref{001})--(\ref{003}) which will improve earlier results obtained recently in~\cite{Cen,SOFMIG,SMH2018,SP}.
Also, an interesting topic is to study optimal control problems for the system (\ref{001})--(\ref{003}) including the necessary conditions of optimality for the control problems, see~\cite{Migorski2020}. 
Finally, it is an extensive and important project to provide the numerical analysis of the system, see~\cite{SHM} and the references therein.


\begin{thebibliography}{99}
	
\baselineskip=13pt

\bibitem{Anh}
N.T.V. Anh, T.D. Ke, 
On the differential variational inequalities of parabolic-parabolic type, 
{\it Acta Appl. Math.} {\bf 176} (2021), 5.

\bibitem{Bartosz}
K. Bartosz, Hemivariational inequalities modeling dynamic contact problems with adhesion, {\it Nonlinear Anal. Theory
Methods Appl.} {\bf 71} (2009), 1747--1762.

	
	
	
	
	
	
	
	
	
	
	
	
\bibitem{CLM}
S. Carl, V.K. Le, D. Motreanu,
{\it Nonsmooth Variational Problems and their Inequalities},
Springer, New York, 2007.

\bibitem{Cen}
J.X. Cen, L. Li, S. Mig\'orski, Van Thien Nguyen, Convergence of a generalized penalty and regularization method for quasi-variational-hemivariational inequalities, {\it Communications in Nonlinear Science and Numerical Simulation} {\bf 103} (2021), 105998. 
	
\bibitem{Clarke}
F.H. Clarke, {\it Optimization and Nonsmooth Analysis}, Wiley, New York, 1983.
	
	
	
	
	
	
	\bibitem{DMP1}
	Z. Denkowski, S. Mig\'orski, N.S. Papageorgiou,
	{\it An Introduction to Non\-li\-near Analysis: Theory}, Kluwer Academic/Plenum Publishers, Boston, Dordrecht, London, New York, 2003.
	
	\bibitem{DMP2}
	Z. Denkowski, S. Mig\'orski, N.S. Papageorgiou,
	{\it An Introduction to Non\-li\-near Analysis: Applications}, Kluwer Academic/Plenum Publishers, Boston, Dordrecht, London, New York, 2003.
	
	
	
	
	
	
	\bibitem{Go11.I}
	D. Goeleven, D. Motreanu, Y. Dumont, M. Rochdi,
	{\em Variational and Hemivariational Inequalities, Theory, Methods and Applications, Volume I: Unilateral Analysis and Unilateral Mechanics}, Kluwer Academic Publishers, Boston,
	Dordrecht, London, 2003.
	
	
	\bibitem{Gwinner}
	J. Gwinner, On a new class of differential variational inequalities and a stability result, 
	{\it Math. Program.} {\bf 139} (2013), 205--221.

\bibitem{HLM2015}
J. Han, Y. Li, S. Mig\'orski, 
Analysis of an adhesive contact problem for viscoelastic materials with long memory, 
{\it Journal of Mathematical Analysis and Applications}  {\bf 427} (2015), 646--668. 
	
	
	
	
	
	\bibitem{HMS2017} 
	W. Han, S. Mig\'orski, M. Sofonea, 
	Analysis of a general dynamic history-dependent variational-hemivariational inequality, 
	{\it Nonlinear Analysis: Real World Applications}
	{\bf 36} (2017), 69--88. 
	
	\bibitem{HS}
	W. Han, M. Sofonea, {\it Quasistatic  Contact Problems in Viscoelasticity and Viscoplasticity}, 
	Studies in Advanced Mathematics {\bf 30}, Americal Mathematical Society, Providence, RI--International Press, Somerville, MA, 2002.
	
	
	
	
	
	
	
	
	\bibitem{KULIG}
	A. Kulig, S. Mig\'orski,
	Solvability and continuous dependence results for second order
	nonlinear inclusion with Volterra-type operator,
	{\it Nonlinear Analysis} {\bf 75} (2012), 4729--4746.
	

\bibitem{Liu2008}
Z.H. Liu, Existence results for quasilinear parabolic hemivariational inequalities, {\it J. Dif\-fe\-rential Equations} {\bf 244} (2008), 1395--1409.

\bibitem{LMZ2017}
Z.H. Liu, S. Mig\'orski, S.D. Zeng, 
Partial differential variational inequalities involving nonlocal boundary conditions in Banach spaces, 
{\it J. Dif\-fe\-rential Equations} {\bf 263} (2017), 3989--4006.
	
\bibitem{Liu1}
Z.H. Liu, S.D. Zeng, D. Motreanu, 
Partial differential hemivariational inequalities, {\it Adv. Nonlinear Anal.} {\bf 7} (2017), 571--586.
	
	

\bibitem{Miettinen}
M. Miettinen, A parabolic hemivariational inequality,  {\it Nonlinear Analysis: Theory, Me\-thods and Applications}
{\bf 26} (1996), 725--734.
	
	\bibitem{MMM}
	S. Mig\'orski, 
	Evolution hemivariational inequality for a class of dynamic viscoelastic nonmonotone 
	frictional contact problems,
	{\it Computers \& Mathematics with Applications} {\bf 52} (2006), 677--698.

\bibitem{Migorski2020}
S. Mig\'orski,
Optimal control of history-dependent evolution inclusions with applications to frictional contact, {\it J.  Optimization Theory and Applications} {\bf 185} (2020), 574--596. 

\bibitem{Migorski2021}
S. Mig\'orski, 
A class of history-dependent systems of evolution inclusions with applications, 
{\it Nonlinear Analysis: Real World Applications}
{\bf 59} (2021), 103246.
	

\bibitem{MHZ2020}
S. Mig\'orski, W. Han, S.D. Zeng, 
A new class of hyperbolic variational-hemivariational inequalities driven by non-linear evolution equations, {\it European J. Appl. Math.} {\bf 32} (2021), 59--88.
	

\bibitem{MO2004}
S. Mig\'orski, A. Ochal, Boundary hemivariational inequality of parabolic type, {\it Nonlinear Analysis: Theory Methods and Appl.} {\bf 57} (2004), 579--596.

\bibitem{MO2008}
S. Mig\'orski, A. Ochal, Dynamic bilateral contact 
problem for viscoelastic piezoelectric materials with adhesion, {\it Nonlinear Analysis: Theory, Methods and Applications} {\bf 69} (2008), 495--509.
	
	\bibitem{MOS13}
	S. Mig\'orski, A. Ochal, M. Sofonea, History-dependent subdifferential inclusions and hemivariational inequalities in contact mechanics, {\it Nonlinear Analysis: Real World Applications} {\bf 12} (2011), 3384--3396.
	
	\bibitem{MOSBOOK}
	S. Mig\'orski, A. Ochal, M. Sofonea,
	{\it Nonlinear Inclusions and Hemivariational Inequalities. Models and Analysis of Contact Problems}, Advances in Mechanics and Mathematics \textbf{26}, Springer, New York, 2013.
	
	\bibitem{MOS18}
	S. Mig\'orski, A. Ochal, M. Sofonea, History-dependent variational-hemi\-va\-ria\-tio\-nal inequalities in contact mechanics, {\it Nonlinear Analysis: Real World Applications} {\bf 22} (2015), 604--618.
	
	\bibitem{AMMA2}
	S. Mig\'orski, A. Ochal, M. Sofonea, Evolutionary inclusions and hemivariational inequalities, 
	Chapter 2 in \emph{Advances in Variational and Hemivariational Inequalities:
		Theory, Numerical Analysis, and Applications}, edited by W. Han, et al., Advances in Mechanics and Mathematics Series, vol. 33 (2015), 39--64, Springer.
	
	
	

\bibitem{MOgorzaly}
S. Mig\'orski, J. Ogorzaly, Dynamic history-dependent variational-hemivariational inequalities with applications to contact mechanics, 
{\it Zeitschrift f\"ur angewandte Mathematik und Physik} {\bf 68} (2017), Article ID.15, 22p.

\bibitem{MigorskiBiao}
S. Mig\'orski, B. Zeng, 
A new class of history--dependent evolutionary variational--he\-mi\-variational inequalities with unilateral constraints, 
{\it Applied Mathematics \& Optimization}
{\bf 84} (2021), 2671--2697.

\bibitem{MZJOGO2018}
S. Mig\'orski, S. Zeng, A class of differential hemivariational inequalities in Banach spaces, 
{\it J. Global Optimization} {\bf 72} (2018), 
761--779.

\bibitem{MShengda2018}
S. Mig\'orski, S.D. Zeng, Hyperbolic hemivariational inequalities controled by evolution equations with application to adhesive contact model, 
{\it Nonlinear Analysis: Real World Applications} 
{\bf 43} (2018), 121--143. 

	
	\bibitem{NP}
	Z. Naniewicz, P.D. Panagiotopoulos,
	{\it Mathematical Theory of Hemivariational Inequalities and Applications}, Marcel Dekker, Inc., New York, Basel, Hong Kong, 1995.
	
	
	
	
	
	
	
	\bibitem{Pang}
	J.S. Pang, D.E. Stewart, Differential variational inequalities, {\it Math. Program.} {\bf 113} (2008), 345--424.
	
	
	
	\bibitem{SST}
	M. Shillor, M. Sofonea, J.J. Telega, {\it Models and Analysis of Quasistatic Contact}, 
	Lect. Notes Phys. {\bf 655}, Springer, Berlin, Heidelberg, 2004.
	
	
	\bibitem{SHM}
	M. Sofonea, W. Han, S. Mig\'orski, 
	Numerical analysis of history-dependent va\-ria\-tio\-nal-hemivariational 
	inequalities with applications to contact problems, 
	{\it European Journal of Applied Mathematics} {\bf 26} (2015), 427--452. 

\bibitem{SHS}
M. Sofonea, W. Han, M. Shillor, 
{\it Analysis and Approximation of Contact Problems with Adhesion or Damage}, Pure and Applied Mathematics
{\bf 276}, Chapman-Hall/CRC Press, New York, 2006.
	
	\bibitem{SM1} 
	M. Sofonea, A. Matei, History-dependent quasivariational
	inequalities arising in Contact Mechanics, {\it European Journal of Applied Mathematics} {\bf 22} (2011), 471--491.
	
	
	\bibitem{SOFMIG}
	M. Sofonea, S. Mig\'orski,
	{\it Variational-Hemivariational Inequalities with Applications},
	Chapman \& Hall/CRC, 
	Boca Raton, 2018.
	
	
	\bibitem{SMH2018}
	M. Sofonea, S. Mig\'orski, W. Han, 
	A penalty method for history-dependent variational-hemivariational inequalities, 
	{\it Computers \& Mathematics with Applications} {\bf 75} (2018), 2561--2573. 
	
	\bibitem{AMMA14}
	M. Sofonea, S. Mig\'orski, A. Ochal, Two history-dependent contact problems, Chapter 14 in
	\emph{Advances in Variational and Hemivariational Inequalities: Theory, Numerical Analysis, and Applications}, 
	edited by W. Han, et al., Advances in Mechanics and Mathematics Series, vol.\ 33 (2015), 355--380, Springer.

	\bibitem{SP} 
	M. Sofonea, F. P\v{a}trulescu, 
	Penalization of history-dependent variational inequalities, {\it European Journal of Applied Mathematics} {\bf 25} (2014), 155--176.
	
	\bibitem{SX} 
	M. Sofonea, Y. Xiao,
	Fully history-dependent quasivariational inequalities in contact mechanics, 
	{\it Applicable Analysis} {\bf 95} (2016), 
	2464--2484.
	
	
\bibitem{ZM2021} 
S.D. Zeng, S. Mig\'orski, 
Dynamic history-dependent hemivariational inequalities controlled by evolution equations with application to contact mechanics,
{\it Journal of Dynamics and Differential Equations}, 2021, in press, 
https://doi.org/10.1007/s10884-021-10088-0.
	
	
\end{thebibliography}
\end{document}